\documentclass{article}

\usepackage[utf8]{inputenc}
\usepackage[a4paper, margin=1.5 in]{geometry}

\usepackage[utf8]{inputenc} 
\usepackage[T1]{fontenc}    
\usepackage{hyperref}       
\usepackage{url}            
\usepackage{booktabs}       
\usepackage{amsfonts}       
\usepackage{nicefrac}       
\usepackage{microtype}      
\usepackage{lipsum}
\usepackage{graphicx}
\usepackage{placeins}
\usepackage{amsmath}
\graphicspath{ {./images/} }
\usepackage[dvipsnames]{xcolor}
\usepackage[ruled,vlined,linesnumbered]{algorithm2e}

\usepackage{tikz}
\usetikzlibrary{automata, matrix}

\usepackage{dcolumn}
\newcolumntype{2}{D{.}{}{2.0}}
\newcommand{\cone}{\mbox{cone}}
\renewcommand{\span}{\mbox{span}}

\newcommand{\rk}{\mbox{rk}_{\R}}
\newcommand{\Rp}{\mathbb{R}_+}
\newcommand{\rkp}{\mbox{rk}_{\Rp}}
\newcommand{\rkb}{\mbox{rk}_\B}

\newcommand{\fr}{\mbox{frk}_\B}
\newcommand{\rkbin}{\mbox{rk}_{\BB}}
\newcommand{\rkz}{\mbox{rk}_{\Z_2}}

\newcommand{\R}{\mathbb{R}}
\newcommand{\B}{\mathcal{B}}
\newcommand{\C}{\mathcal{C}}

\newcommand{\Z}{\mathbb{Z}}
\newcommand{\BB}{\mathbb{B}}
\newcommand{\F}{\mathbb{F}}
\newcommand{\SR}{\mathcal{S}}

\newcommand{\TEXTFONTT}[1]{\mathcal{#1}}

\newcommand{\G}{\TEXTFONTT{G}}

\DeclareMathOperator*{\argmin}{arg\,min}

\providecommand{\keywords}[1]
{
  \small	
  \textbf{\textit{Keywords---}} #1
}
\usepackage{amsthm}
\newtheorem{theorem}{Theorem}
\newtheorem{corollary}{Corollary}
\newtheorem{proposition}{Proposition}

\theoremstyle{definition}
\newtheorem{definition}{Definition}

\newtheorem{example}{Example}

\usepackage{color,soul}

\newtheorem{innercustomthm}{Theorem}
\newenvironment{customthm}[1]
  {\renewcommand\theinnercustomthm{#1}\innercustomthm}
  {\endinnercustomthm}
  
\newtheorem{innercustomthmm}{Proposition}
\newenvironment{customthm2}[1]
  {\renewcommand\theinnercustomthmm{#1}\innercustomthmm}
  {\endinnercustomthm}

\title{Factorization of Binary Matrices: Rank Relations, Uniqueness and Model Selection of Boolean Decomposition}
\author{Derek DeSantis\footnote{Theoretical Division -  CNLS, Los Alamos National Laboratory, USA ddesantis@lanl.gov}, Erik Skau\footnote{CCS Division, Los Alamos National Laboratory, USA ewskau@lanl.gov}, Duc P. Truong\footnote{CCS Division, Los Alamos National Laboratory, USA dptruong@lanl.gov}, Boian Alexandrov\footnote{Theoretical Division, Los Alamos National Laboratory, USA boian@lanl.gov}}

\date{\today}

\begin{document}

\maketitle

\begin{abstract}
  The application of binary matrices are numerous. Representing a matrix as a mixture of a small collection of latent vectors via low-rank decomposition is often seen as an advantageous method to interpret and analyze data. In this work, we examine the factorizations of binary matrices using standard arithmetic (real and nonnegative) and logical operations (Boolean and $\Z_2$). We examine the relationships between the different ranks, and discuss when factorization is unique. In particular, we characterize when a Boolean factorization $X = W \land H$ has a unique $W$, a unique $H$ (for a fixed $W$), and when both $W$ and $H$ are unique, given a rank constraint. We introduce a method for robust Boolean model selection, called BMF$k$, and show on numerical examples that BMF$k$ not only accurately determines the correct number of Boolean latent features but reconstruct the pre-determined factors accurately. 
\end{abstract}


\section{Introduction}

The rank of a real valued $N \times M$ matrix $X \in \R^{N,M}$ is the dimensionality of the vector space spanned by its columns. If the rank of $X$ is $R$, then the matrix $X$ can be written as a product of two matrices $X = WH$ where $W \in \R^{N,R}$ and $H \in \R^{R,M}$. In this fashion, the columns of $X$ can be seen as a mixing of latent features $w_1, \dots, w_R$, the columns of $W$, according to their weights, the columns of $H$. By imposing various constraints, one obtains different factorizations. Perhaps the  best-known factorization is the Singular Value Decomposition \cite{stewart1993early}, where the factor matrices are restricted to be orthogonal. When the elements of the matrix $X$ admit specific properties, these properties often suggest appropriate constraints on the factor matrices. 
For example, when $X$ is nonnegative, it is natural to impose a nonnegative constraint, and the nonnegative rank is defined as the smallest number $R$ such that $X=WH$ for nonnegative matrices $W \in \R_+^{N,R}$ and $H \in \R_+^{R,M}$ \cite{cohen1993nonnegative}. Similarly, binary rank is the smallest number $R$ for which  a binary matrix can be decomposed into a product of binary matrices \cite{zhang2007binary}. Examples of applications of binary decompositions include social networks, market-basket data, DNA transcription profiles, and many others \cite{li2005general}.

Instead of applying constraints, one can choose to change the underlying arithmetic to arrive at new types of decompositions and ranks. If one changes the arithmetic operations of ``plus'' and ``times'' to logical operations of ``or'' with ``and'', this results in Boolean rank and Boolean matrix factorization (BMF) \cite{miettinen2014mdl4bmf}. Boolean decompositions solve the tiling problem \cite{geerts2004tiling} that determines how to cover the 1's of a binary matrix by a minimum number of subitems, which is equivalent to the NP-hard bi-clique cover problem \cite{doherty1999biclique,orlin1977contentment}.  One can also choose the logical operations of ``xor'' with ``and'', which results in the Galois field decompositions \cite{gutch2012ica,yeredor2011independent}. 

In all the cases discussed above, matrix factorization allows one to learn latent factors from a complex data subject to various constraints and relations between the elements of the data. Many applications in machine learning and data mining, e.g., document classification, recommendation systems, community detection, cryptography and others, involve data with binary values $\{0, 1\}$ \cite{miron2021boolean}. As such,  matrices consisting of 0's and 1's arise in various domains of applicability, and discovering their latent structure is critical for doing any fundamental analysis.

In this work, we present the mathematical theory allowing for comparisons between 1) real, nonnegative, Boolean, binary and $\Z_2$ ranks and 2) the uniqueness of the corresponding factorizations in these  different contexts. The primary focus is to build a framework on which to extract not only the correct number of hidden features, but the correct ones when uniqueness is present. Unfortunately, exact factorizations aren't always possible and even single bit flips can destroy the uniqueness (see Example \ref{ex: destroy unique}).  Practically this leads one to naturally desire a model selection algorithm which 1) discovers a stable rank such that 2) uniqueness is approximately recovered.

In Section 2, we define the notations and standardize the common definitions used throughout the text.  In Section 3, we provide the relationships between the five ranks and examples which demonstrate that no other relationships can exist. In Section 4, we move on to discuss the uniqueness of these types of decompositions, which contains our main results. We begin by redefining a familiar geometric framework for nonnegative factorizations to also include real, $\Z_2$, and Boolean factorizations. We remark that while the real and $\Z_2$ factorizations are always (essentially) unique, nonnegative and Boolean need not to be. Indeed, given a Boolean factorizations, $X = W \cdot H$, of a $0-1$ matrix $X$, we note that uniqueness can be achieved in the patterns $W$ without unique feature weights, $H$, and vice versa.  Hence, we investigate criterion for which $W$, $H$ or $W$ and $H$ are unique. For non-negative factorizations, these results either follow naturally from linear algebra or are well known in the literature.  However in Boolean, these results are new.  

Proposition \ref{prop: unique W} states that $W$ is unique if and only if there exists a unique cone, an additive set analogous to a subspace, which contains the data $X$. We define a property called \textit{freeness} analogous to linear independence, and show in Theorem \ref{thm: unique H if and only if} that the columns of $X$ satisfy a freeness property if and only if the feature weights $H$ are unique. Freeness leads naturally to the definition of a free rank for Boolean matrices. We then show that if the free rank is equal to the Boolean rank, that uniqueness is guaranteed. 

In Section 5, we introduce a method for robust Boolean model selection - BMF$k$. This is a Boolean analog of the nonnegative matrix model selection algorithm NMF$k$ \cite{alexandrov2013signatures,vangara2021finding}. We compare BMF$k$ to NMF$k$ and demonstrate that for a matrix $X$ with a unique Boolean factorization and $\rkb(X) \neq \rkp(X)$, BMF$k$ discovers the correct hidden patterns while NMF$k$ does not. Using theory developed in Section 4, we construct a set of Boolean matrices with unique Boolean factorizations and show that not only does BMF$k$ correctly identify the latent dimension, but it also accurately extract the predetermined latent features. Adding noise to this set, we find that BMF$k$ finds features that are highly correlated with the unique  ``true'' features, which is important for a practical Boolean model selection and extraction technique . 



\section{Definitions}

Throughout, we let $\BB := \{0,1\}$.  When dealing with a $N \times M$ binary matrix $X \in \BB^{M , N}$ one can consider different decompositions $X = W  H$ where the pattern matrix $W$ and weight matrix  $H$ either belong to different sets, such as the reals $\R$, the non-negatives $\R_+$, or the binary set $\BB$. Or alternatively one may want to consider decompositions employing different algebraic operations in the matrix multiplication of $W$ and $H$. Three natural "addition" operators are the addition, $+$, exclusive or $\oplus$, and logical or $\lor$:
\[
\renewcommand\arraystretch{1.3}
\setlength\doublerulesep{0pt}
\begin{tabular}{r||*{2}{2|}}
$+$ & 0 & 1  \\
\hline\hline
0 & 0 & 1  \\ 
\hline
1 & 1 & 2  \\ 
\hline
\end{tabular}
\hspace{.5 in}
\begin{tabular}{r||*{2}{2|}}
$\oplus$ & 0 & 1  \\
\hline\hline
0 & 0 & 1  \\ 
\hline
1 & 1 & 0  \\ 
\hline
\end{tabular}
\hspace{.5 in}
\renewcommand\arraystretch{1.3}
\setlength\doublerulesep{0pt}
\begin{tabular}{r||*{2}{2|}}
$\lor$ & 0 & 1  \\
\hline\hline
0 & 0 & 1  \\ 
\hline
1 & 1 & 1  \\ 
\hline
\end{tabular}.
\]
These addition operators are typically paired with multiplication $\times$, logical and $\otimes$, and logical and $\land$ respectively:
\[
\renewcommand\arraystretch{1.3}
\setlength\doublerulesep{0pt}
\begin{tabular}{r||*{2}{2|}}
$\times$ & 0 & 1  \\
\hline\hline
0 & 0 & 0  \\ 
\hline
1 & 0 & 1  \\ 
\hline
\end{tabular}
\hspace{.5 in}
\begin{tabular}{r||*{2}{2|}}
$\otimes$ & 0 & 1  \\
\hline\hline
0 & 0 & 0  \\ 
\hline
1 & 0 & 1  \\ 
\hline
\end{tabular}
\hspace{.5 in}
\renewcommand\arraystretch{1.3}
\setlength\doublerulesep{0pt}
\begin{tabular}{r||*{2}{2|}}
$\land$ & 0 & 1  \\
\hline\hline
0 & 0 & 0  \\ 
\hline
1 & 0 & 1  \\ 
\hline
\end{tabular},
\]
which are all identical on the binary set $\BB$.  Consequently, there are several different types of arithmetic structures one can consider on binary vectors and matrices. Throughout, we adopt the following notation:

\begin{itemize}
    \item The \textit{real numbers} $\R$ is equipped with the the operations $(+, \times)$
    \item The \textit{nonnegative real numbers} $\R_+ = \{x \in \R: x\geq 0\}$ is equipped with the the operations $(+, \times)$
    \item The \textit{Galois field} $\Z_2 = \BB$ is equipped with the the operations $(\oplus, \otimes)$
    \item The \textit{Booleans} $\B = \BB$ is equipped with the the operations $(\lor, \land)$
\end{itemize}
In what follows, the notation $\R,\R_+,\Z_2$ and $\B$ will always imply usage of the associated operations. We will also make limited use of the restriction of $\R_+$ to $\BB$, though this is not closed under addition. 

 Each of the three pairs of arithmetic operations $(+, \times), (\oplus, \otimes)$ and $(\lor, \land)$ naturally defines matrix multiplication. Given two matrices $W \in \BB^{N,R}$ and $H \in \BB^{R,M}$, we use $WH$, $W \otimes H$, and $W \land H$ to denote the real, $\Z_2$, and Boolean matrix multiplications respectively. When referring to more than one potential operation, we will utilize the notation $W \cdot H$. These different constraints and algebraic operators lead to five natural factorizations of $0-1$ matrices:

\begin{definition}
Let $X \in \BB^{N, M}$.
\begin{itemize}
    \item The \textit{real rank} is $\rk(X) := \min \{R | X = WH, W \in \R^{N,R}, H \in \R^{R,M}\}$.
    \item The \textit{nonnegative rank} is $\rkp(X) := \min \{R | X = WH, W \in \R_+^{N,R}, H \in \R_+^{R,M}\}$.
    \item The \textit{binary rank} is $\rkbin(X) := \min \{R | X = WH, W \in \BB^{N,R}, H \in \BB^{R,M}\}$.
    \item The \textit{$\Z_2$ rank} is $\rkz(X) := \min \{R | X = W \otimes H, W \in \Z_2^{N,R}, H \in \Z_2^{R,M}\}$.
    \item The \textit{boolean rank} is $\rkb(X) := \min \{R | X = W \land H, W \in \B^{N,R}, H \in \B^{R,M}\}$.
\end{itemize}
For each defined rank \footnote{The definitions of rank given here correspond to what is often referred to as the \textit{factor} or \textit{Schein} rank.}, a factorization is said to be a \emph{rank factorization} or \emph{rank revealing} if $W$ and $H$ correspond to a minimum $R$ factorization. 
\end{definition}

Recall that a field is a set wherein addition, subtraction, multiplication and division are well defined, and a semiring is one where only addition and multipication are well defined. We note that $\R$ and $\Z_2$ are fields, while $\R_+$ and $\B$ are semirings which are not fields.  Hence, $\R^N$ and $\Z_2^N$ are vector spaces over their respective fields, and $\R_+^N$ and $\B^N$ are semimodules over their respective semirings \cite{gondran2008graphs, roman2005advanced}.  Therefore, they are afforded their respective properties, which we will infrequently leverage. Binary factorizations are not a closed arithmetic, and therefore can only be viewed as the restriction of $\R_+$ to $\{0,1\}$. 

Our goal is to discuss rank relations as well as uniqueness results for $0-1$ matrices. For historical reasons, the literature has generally compared factorizations against the reals. However these comparisons are sometimes artificial, so we take the more mathematically consistent approach of considering the two cases of field factorizations and semiring factorizations. 

\section{Rank relations}

In this section, we briefly discuss how the five different ranks for $0-1$ matrices are related.  In some cases, two ranks are incomparable in the sense that one is not guaranteed to be larger than the other. However, there are four guaranteed relations between these five ranks which are summarized in Figure \ref{fig:rank inequalities}:

\begin{figure}[ht!]
\begin{center}
\begin{tikzpicture}[
        > = stealth, 
        shorten > = 1pt, 
        auto,
        node distance = 3cm, 
        semithick 
    ]
    every state/.style={%
        draw = black,
        thick,
        fill = white,
        minimum size = 1mm
    }
    \matrix[%
        matrix of math nodes,
        column sep = 2cm,
        row sep = 1cm,
        inner sep = 0pt,
        nodes={state}
        ] (m) {%
        \rkz & \rk & \rkp & \rkbin \\
                &        &  \rkb&   \\
        };
        \path[->] (m-1-1) edge node {$\leq$} (m-1-2)
                  (m-1-2) edge node {$\leq$} (m-1-3)
                  (m-1-3) edge node {$\leq$} (m-1-4)
                  (m-2-3) edge node [pos=0.25,above right, rotate=90] {$\leq$} (m-1-3);
\end{tikzpicture}
\end{center}
    \caption{Flowchart of the five rank inequalities.}
    \label{fig:rank inequalities}
\end{figure}
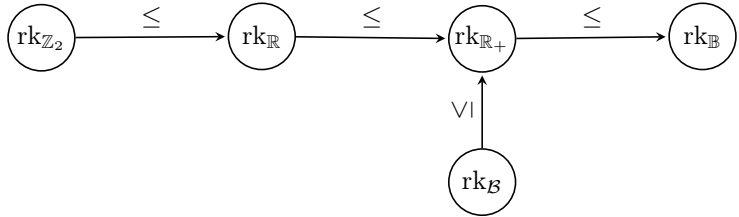

Since nonnegative factorization is a restriction of real, and binary is a restriction of nonnegative, $\rk(X) \leq \rkp(X) \leq \rkbin(X)$ follow by definition. The relationship $\rkb(X) \leq \rkp(X)$ was previously known, as mentioned in \cite {miron2021boolean}, and $\rkz(X) \leq \rk(X)$ appears to be a new result, which is easy to prove by inspecting determinants.  Since we could not find a detailed proof of either of these results, we state these results and provide the proofs in the appendix for completeness.

To prove $\rkb(X) \leq \rkp(X)$ we first define some notations and relations. Given a nonnegative matrix $X = (X_{i,j})$, we can form a binary matrix $\hat{X}$, via: $\hat{X}_{i,j} = 1$ if $X_{i,j} > 0$, and $\hat{X}_{i,j} = 0$ otherwise. This technique is often referred to as \textit{thresholding} \cite{zhang2007binary}.  The following proposition states that thresholding a nonnegative decomposition reveals a Boolean decomposition:

\begin{proposition}
\label{prop: NMF to Boolean}
Let  $X \in \BB^{N,M}$.  If $X = WH$ is a nonnegative factorization of $X$, then $\hat{W} \land \hat{H}$ is a Boolean factorization of $X$. 
\end{proposition}

Proposition \ref{prop: NMF to Boolean} states that the process of converting the nonnegative factors $W,H$ of a binary matrix $X$ into Boolean matrices $\hat{H}, \hat{W}$, via thresholding produces a Boolean factorization. However, we remark that this factorization may not be a minimal rank Boolean factorization.  This is the contents of the following theorem:

\begin{theorem}
Let $X \in \BB^{N,M}$. Then $\rkb(X) \leq \rkp(X)$.
\end{theorem}

\begin{proof}
By Proposition \ref{prop: NMF to Boolean}, each nonnegative factorization yields a Boolean factorization.  In particular if $K = \rkp(X)$, and $X = WH$ is a nonnegative rank $K$ factorization of $X$, then $\hat{W} \land \hat{H}$ is a rank $K$ Boolean factorization of $X$.  Thus, $\rkb(X) \leq \rkp(X)$.
\end{proof}

Next, we state the connection between the $\Z_2$ rank and the real rank.  Its proof can be found in the Appendix.

\begin{theorem}
\label{thm: z2 rank}
Let $X \in \BB^{N,M}$. Then $\rkz(X) \leq \rk(X)$.
\end{theorem}

Given $X \in \BB^{N,M}$, two entries $X_{i,j}$ and $X_{k,l}$ are called \textit{isolated ones} if $X_{i,j} = X_{k,l} =1 $ and $X_{i,l} X_{k,j}=0$ \cite{beasley2012isolation}. For example, if 
$$X = \begin{bmatrix}
1 & 0 & 1\\
0 & 1 & 0\\
0 & 0 & 1
\end{bmatrix},$$
then $\{X_{1,1}, X_{2,2},X_{3,3}\}$ are isolated ones. The \textit{isolation number} $\iota(X)$ is the maximum size of all sets of isolated ones \cite{beasley2012isolation}. The following useful proposition provides a lower bound on some the ranks of a matrix in terms of the isolation number.

\begin{proposition}[\cite{beasley2012isolation}]
\label{prop: pairwise ind}
Let $X \in \BB^{N,M}$.  Then $$\iota(X) \leq \rkb(X) \leq \rkp(X) \leq \rkbin(X).$$
\end{proposition}

\subsection{Rank examples}

We will show that the inequalities in Figure \ref{fig:rank inequalities} are strict, and that no other possible rank comparisions hold.  For this, we will utilize the following series of examples:

\begin{example}
\label{ex: A}
Let 
$$A = \begin{bmatrix}
0 & 1 & 1\\
1 & 0 & 1\\
1 & 1 & 0
\end{bmatrix}.$$
It is easy to see that $\rkz(A) =2$ and $\rk(A) = 3$. On the other hand, $A$ has three isolated ones, $\{A_{1,2}, A_{2,3}, A_{3,1}\}$, so by Proposition \ref{prop: pairwise ind} the Boolean rank is $3$. 
\end{example}

\begin{example}
\label{ex: B}
Let
$$B = \begin{bmatrix}
0 & 1 & 1\\
1 & 0 & 1\\
1 & 1 & 1
\end{bmatrix}.$$
Then each column of $B$ are $\Z_2$ independent.  Moreover, since none of the columns sum to another, we have that $\rkz(B) = 3$.   However, $(0,1,1) \vee (1,0,1) = (1,1,1)$.  Thus, $\rkb(B) \leq 2$. Since $B_{1,2}$ and $B_{2,1}$ are isolated ones, we have that $\rkb(B) = 2$.  Furthermore, $\rkp(B) = \rk(B) = 3$. 
\end{example}

\begin{example}[\cite{cohen1993nonnegative}]
\label{ex: C}
Let 
$$C = \begin{bmatrix}
1 & 1 & 0 & 0\\
1 & 0 & 1 & 0\\
0 & 1 & 0 & 1\\
0 & 0 & 1 & 1
\end{bmatrix}$$
Then the entries $\{ C_{2,1}, C_{1,2}, C_{3,4}, C_{4,3}\}$ are isolated ones.  Hence by Proposition \ref{prop: pairwise ind} $\rkb(C) = 4$, so that $\rkp(C) = 4$, while, $\rk(C) = 3$. 
\end{example}

\begin{example}
\label{ex: D}
Let 
$$D = \begin{bmatrix}
1 & 1 & 1 & 1 & 1 & 1\\
1 & 1 & 0 & 0 & 1 & 0\\
1 & 0 & 1 & 0 & 0 & 1\\
0 & 1 & 0 & 1 & 0 & 1\\
0 & 0 & 1 & 1 & 1 & 0
\end{bmatrix}$$
then $4 =\rkp(D) < \rkbin(D) = 5$.  To see this, note that $D$ contains the matrix $C$ from Example~\ref{ex: C} as a submatrix in the bottom left corner.  Thus, $\rkp(D) \geq 4$. Furthermore, one has the decomposition
\[
D =
\begin{bmatrix}
0.5 & 0.5 & 0.5 & 0.5 \\
1 & 0 & 0 & 0\\
0 & 1 & 0 & 0\\
0 & 0 & 1 & 0\\
0 & 0 & 0 & 1\\
\end{bmatrix}
\begin{bmatrix}
1 & 1 & 0  & 0  & 1  & 0\\
1 & 0 & 1  & 0  & 0 & 1\\
0 & 1  & 0  & 1  & 0  & 1\\
0 & 0  & 1  & 1 & 1 & 0
\end{bmatrix}.
\]
Thus $\rkp(D) \leq 4$.   One can also observe that $\rk(D) = 4$, so that Theorem 6 of \cite{gillis2012sparse} may be applied to show the above decomposition is unique. Since binary factorizations are a restriction of nonnegative, one concludes that $\rkbin(D)> 4$, and therefore $\rkbin(D) = 5$. 
\end{example}

Table \ref{table: ranks} summarizes the relevant information in Examples \ref{ex: A}-\ref{ex: D}. Example matrices $A$, $B$, $C$ and $D$ show that the inequalities in Figure \ref{fig:rank inequalities} are strict. Moreover, example matrices $A$ and $B$ show that $\rkz$ and $\rkb$ are not comparable, while example matrices $B$ and $C$ show that $\rk$ and $\rkb$ are not comparable. 

\begin{table}[h!]
    \centering
    \begin{tabular}{|c||c|c|c|c|c|}
        \hline
        Example &  $\rk$ & $\rkz$ & $\rkp$ & $\rkb$ & $\rkbin$\\
        \hline
        \hline
        $A$ &  $3$ & $2$ & - & $3$ & -\\
        $B$ &  $3$ & $3$ & $3$ & $2$ & -\\
        $C$ &  $3$ & - & $4$ & $4$ & -\\
        $D$ &  - & - & $4$ & - & $5$\\
        \hline
    \end{tabular}
    \caption{Ranks from Examples \ref{ex: A}-\ref{ex: D}.}
    \label{table: ranks}
\end{table}

\section{Uniqueness of rank factorizations and geometry}

Let $X \in \BB^{N,M}$, and $X=W \cdot H$ be any of the five rank factorizations.  If $P$ is an $R \times R$ permutation matrix, then $W' = W \cdot P$ and $H' = P^{-1} \cdot H$ is another rank factorization for $X$ of the same type. More generally, one could use any invertible matrix whose inverse is within the algebra of interest:

\begin{definition}
Given a semiring $\SR$ (which may be a field), a matrix $U \in \SR^{R,R}$ is said to be \textit{invertible} if there exists another matrix $U^{-1} \in \SR^{R,R}$ such that $U U^{-1} = U^{-1} U = I$, the identity matrix.
\end{definition}

If $X \in \BB^{N,M}$, and $X=W \cdot H$ is a rank factorization, and $U$ is invertible over the respective algebra, then $W' = W \cdot U$ and $H' = U^{-1} \cdot H$ is another rank factorization for $X$ of the same type.  We note that if $U \in \R_+^{R,R}$ is an invertible matrix such that $U^{-1} \in \R_+^{R,R}$, then $U$ must be of the form $DP$ for some positive diagonal matrix $D$ and permutation matrix $P$ \cite{brown1968invertibly}. Similarly, if $U \in \B^{R,R}$ is invertible, then it must be the case that $U = P$ for some permutation matrix \cite{ledley1965inverse}.  In either case, multiplication by an invertible matrix encodes a trivial change.  For this reason, we make the following definition for uniqueness:

\begin{definition}
Let $X \in \BB^{N,M}$.  A factorization $X = W \cdot H$ is said to be \textit{unique} if given another factorization $X=W' \cdot H'$, one has that $W' = W \cdot U$ and $H' = U^{-1} \cdot H$ for some matrix $U$ that is invertible over the respective algebra. 
\end{definition}

It is easy to see that rank factorizations over fields are always unique. Recall that given a field $\F$ and a positive integer $N$, the space $\F^N$ defines a vector space with scalars drawn from $\F$. Any linear map $\F^N \rightarrow \F^M$ defines an $N \times M$ matrix with entries from $\F$, and vice versa. The appropriate geometric notion to describe rank factorizations of matrices over a field is a subspace.  Considering the matrix $X$ as a linear map from $\F^N \rightarrow \F^M$, the rank $\mbox{rk}_{\F}(X)$ is equal to the maximal number of $\F-$linearly independent vectors in the range of $X$, denoted $\mbox{rng}(X)$. Therefore if $R = \mbox{rk}_{\F}(X)$, each factorization $X = W \cdot H$ for $W \in \F^{N,R}$, $H \in \F^{R,M}$ is in 1-1 correspondence with a choice $w_1, \dots w_R$ of basis for $\mbox{rng}(X)$. Hence, any two rank factorizations differ by a change of basis, which is encoded in an invertible matrix.

The story for semiring factorizations is similar, but comes with some caveats.   We will show that finding a rank factorization in these cases is equivalent to finding a minimal type of generating set of a \textit{cone} rather than a subspace.  While there is a unique minimal subspace that contains the data in the field case (namely the range of the matrix), there may not be a unique minimal cone that contains the data in the semiring cases. We begin with the definition of a cone:



\begin{definition}
Let $\SR$ be a semiring ($\SR  = (\R_+, +, \times)$ or $(\B,\lor, \land)$) and $N \in \mathbb{N}$. 
A \textit{cone} is a subset $\C \subset \SR^N$ that is closed under addition from $\SR^N$ and scalar multiplication from $\SR$ \footnote{Those familiar with semi-ring theory might note that the definition of a cone is the same as an $\SR$ sub-semimodule over the semiring $\SR^N$.  We have chosen to use the term cone rather than sub-semimodule due to its ubiquity in nonnegative matrix factorizations.}. Given a collection $\{W_i\}_{i=1}^R \subset \SR^N$, we define the \textit{span} of $\{W_i\}_{i=1}^R$ as
\[
\mbox{span} \{W_i\}_{i=1}^R := \left \{ \sum_{i=1}^R h_i \cdot W_i: W_i \in \SR^N, h_i \in \SR \right\}
\]
where scalar multiplication and addition utilize the element-wise operations from $\SR$.  Given a matrix $W \in \SR^{N,R}$, we define the \textit{cone of the matrix $W$} to be
\[
\mbox{cone}(W) = \{W \cdot h: h \in \SR^R\} \subset \SR^N.
\]
A cone $\C$ is said to have a \textit{generating set} $G \subset \C \subset \SR^N$ if $\mbox{span} (G) = \C$. The \textit{order} of the cone $\C \subset \SR^N$, denoted $\mathcal{O}(\C)$ is the size of a minimal generating set. 
A cone is \textit{finitely generated} if $\mathcal{O}(\C) < \infty$. 
\end{definition}

Every finitely generated  cone $\C \subset \SR^N$ is easily seen to be the  span of some collection of vectors and vice versa.  Furthermore, $\mbox{cone}(W)$ is a finitely generated cone, and every finitely generated cone is $\mbox{cone}(W)$ for some matrix $W$. We will only be interested in finitely generated cones, so we drop the text \textit{finitely generated} in what follows. We also remark that for $\SR = (\R_+, +, \times)$, a cone $\C$ is equivalently the intersection of half spaces, though we don't use this in our analysis. 



A cone $\C \subset \SR^N$ may have many different generating sets.  For example, consider the order 2 cone $\C \subset \B^2$ given by 
\[
\C = \mbox{cone}\left(
\begin{pmatrix}
1 & 0 \\
0 & 1 
\end{pmatrix}
\right)
=
\left\{
\begin{pmatrix} 
0\\
0
\end{pmatrix},
\begin{pmatrix}
1\\
0
\end{pmatrix},
\begin{pmatrix}
0\\
1
\end{pmatrix},
\begin{pmatrix}
1\\
1
\end{pmatrix}
\right\} =
\mbox{cone}\left(
\begin{pmatrix}
1 & 0 & 1 \\
0 & 1 & 1
\end{pmatrix}
\right).
\]

Any compact convex subset of $\R_+^N$ is equal to the convex hull of the extreme rays by the Krien-Milman theorem.  Consequently a minimal generating set for a cone $\C \subset \R_+^N$ can be uniquely described up to positive scalar multiplication by the extreme rays. It turns out that for any cone $\C$ over $\B$, the minimal generating set is also unique.  Thus in our two semiring cases of interest, $\R_+$ and $\B$, we can unambiguously talk about \textit{the} minimal generating set for the cone $\C$.  See the appendix for details. 

As mentioned above, in fields each factorization $X = W \cdot H$ for $W \in \F^{N,R}$, $H \in \F^{R,M}$ is in 1-1 correspondence with a basis for $\mbox{rng}(X)$.  The same is true for matrix factorizations over a semirings $\R_+$ and $\B$.  The following proposition is a consequence of the definition:

\begin{proposition}
\label{cone iff factorization}
Let $X \in \SR^{N,M}$ for $\SR = (\R_+, +, \times)$ or $(\B,\lor, \land)$.  Then $X$ has a factorization $X = W \cdot H$ for some $W \in \SR^{N,R}$ and $H \in\SR^{R,M}$ if and only if $\cone(X) \subset \cone(W)$.  Moreover, 
\[
\mbox{rk}_{\SR}(X) = \min \{\mathcal{O}(\C): \C \subset \SR^N \mbox{ is a cone with } \cone(X) \subset \C\}.
\]
\end{proposition}




Proposition \ref{cone iff factorization} states that finding a rank factorization $X=W \cdot H$ is equivalent to finding a cone $\C$ of minimal order that contains the data $\mbox{cone}(X)$.  In general, there can be more than one such minimal cone $\C$, as the following example demonstrates:

\begin{example}
\label{ex: bad matrix}
Consider nonnegative factorizations of the matrix
$$X = \begin{bmatrix}
1 & 1 & 0 & 0\\
1 & 0 & 1 & 0\\
0 & 1 & 0 & 1\\
0 & 0 & 1 & 1
\end{bmatrix}.$$
It is well known that $X$ has two distinct nonnegative rank 4 NMF factorizations: $X = IX$ and $X= XI$, where $I$ is the $4 \times 4$ identity.  As such, there are two separate sets of extreme rays that can be used to create the data $X$. From Example \ref{ex: C} we know that $\rkb(X) = 4$, so this fails to have a unique Boolean rank factorization as well.
\end{example}

Another issue that can arise in the semi-ring cases is non-uniqueness of the feature weights, $H$:
\begin{example}
\label{ex: unique W, not unique H}
Consider Boolean decompositions of the Boolean matrix
\[
X = \begin{bmatrix}
1 & 1 & 0 & 1\\
1 & 0 & 1 & 1\\
0 & 0 & 1 & 1\\
0 & 1 & 1 & 1
\end{bmatrix}.
\]
If $x_i$ is the $i$'th column of $X$, then $W = [x_1 \ x_2 \  x_3]$ is the unique minimal generating set for $\mbox{cone}(X)$.  However, since $x_1 \vee x_3 = x_2 \vee x_3 = x_4$, we see that the weights $H$ cannot be unique. 
\end{example}

While field factorizations are always unique, the above examples demonstrate that semiring factorizations can fail to be unique because either $W$, $H$ or both $W$ and $H$ can be not unique.  We now investigate criterion that guarantee uniqueness of one of the factors, or both of the factors. 

\subsection{Uniqueness of the patterns $W$}

We begin by discussing how to achieve uniquness of the factor $W$.  From Proposition \ref{cone iff factorization}, one immediately achieves the following uniqueness statement for the patterns, $W$ in terms of geometry of the cones:

\begin{proposition}
\label{prop: unique W}
Let $X \in \SR^{N, M}$ have $\mbox{rk}_{\SR}(X) = R$ for $\SR = \R_+$ or $\B$. Then $X=W \cdot H$ for $W \in \SR^{N,R}$, $H \in \SR^{R,M}$ has a unique $W$ if and only if there exists a unique order $R$ cone $\C \subset \SR^N$ such that $\cone(X) \subset \C$.
\end{proposition}



While Proposition \ref{prop: unique W} characterizes the uniqueness of cones, it isn't always the most practical to implement. Given a $W$, it would be nice to know if it is indeed unique. In \cite{laurberg2008theorems}, the authors defined a condition called boundary close which was a necessary condition for uniqueness of $W$ in NMF. With very mild changes, these conditions also hold in Boolean case.

\begin{definition}
 A subset $T \subset \BB^R$ is \textit{boundary close} if for each $i,j \in \{1,\dots, R\}$ with $i \neq j$, there exists a $t \in T$ such that $t_i = 1$ and $t_j = 0$.
\end{definition}

Notice that the row vectors for $W$ in Example \ref{ex: unique W, not unique H} are boundary close. This happens to always be the case when there is a unique minimal cone containing the data $X$:

\begin{proposition}
Let $X = W \cdot H$ be a nonnegative or Boolean factorization. Then if $\mbox{cone}(W)$ is the unique order $R$ cone which contains $\cone(X)$, then the set of row vectors for $W$ are boundary close.
\end{proposition}

The proof strategy for NMF implemented in \cite{laurberg2008theorems} holds for Boolean with very mild changes, so we will not reproduce it here. While boundary close is necessary for uniqueness, it is easily seen to not be sufficient, e.g., see Example \ref{ex: bad matrix}.

\subsection{Uniqueness of the feature weights H}

Suppose that $\cone(X) \subset \cone(W)$ for some matrix $W$.  Then by Proposition \ref{cone iff factorization}, $X = W \cdot H$ for some feature weights $H$.  This  matrix may not be unique even if $W$ is. We now discuss how to achieve uniqueness of the weights matrix $H$ given a fixed pattern matrix $W$.  We will show that when each column of $X$ satisfies a particular property in terms of $W$, then $H$ must be unique. For this we investigate NMF and Boolean separately, as the conditions are stated slightly different.

\subsubsection{NMF}
As $\R_+^N \subset \R^N$, one can simply utilize linear independence to codify the uniqueness of a weight matrix $H$. Recall that the solution space of the real linear system $x = Wh$ is given by $h = W^\dagger x + q$, where $W^\dagger$ is the Moore-Penrose inverse of $W$ and $q$ is any element in the null space of $W$. The set  $\{W^\dagger x +q: q \in \operatorname{null}(W) \}$ intersects the nonnegative orthant at a singleton if and only if $x=Wh$ has a unique nonnegative solution. This is summarized in the following proposition:

\begin{proposition}
\label{prop: NMF unique H if and only if}
Suppose $X \in \R_+^{N, M}$, $W \in \R_+^{N,R}$ and $\cone(X) \subset \cone(W)$. Then $X = WH$ has a unique $H \in \R_+^{R,M}$ for the fixed $W$ if and only if $\{ W^\dagger x_i + q : q \in \operatorname{null}(W) \}  \cap \R_+^R$ is a singleton for each column $x_i$.
\end{proposition}

If $\cone(X) \subset \cone(W)$ and $W$ has full rank, then the $\operatorname{null}(W) = 0$.  Thus, the solution space is a singleton which is nonnegative by assumption.  We record this in the following corollary: 

\begin{corollary}
\label{corollary: NMF unique H if}
Suppose $X \in \R_+^{N, M}$, $W \in \R_+^{N,R}$ and $\cone(X) \subset \cone(W)$. If $\rk(W) = R$, then $X = W H$ has a unique $H$  for the fixed $W$.
\end{corollary}

The converse of Corollary \ref{corollary: NMF unique H if} does not hold.  Indeed, if $X$ is as in Example \ref{ex: bad matrix} with $W = X$, then $H = I$ is unique.  However, $\rk(W) = 3 < \rkp(W) = 4$. We remark that in the case when $\rk(X) = \rkp(X)$, uniqueness of $W$ guarantees uniqueness of $H$.  This is the contents of Theorem \ref{thm: NMF unique} below.


\subsubsection{Boolean}
$\B^N$ does not naturally reside inside a field $\F^N$ in the same way that $\R_+$ sits inside of $\R$ \footnote{If it did, then since addition over a field forms an invertible group, each vector must have an additive inverse.  But each Boolean vector is idempotent, and the only invertible idempotent in a group zero element.}.  Consequently, concepts such as linear independence no longer hold. While there is a concept of \textit{independence} for a generating set of a semiring, we will require a slightly stronger criterion to investigate the uniqueness.

\begin{definition}
A collection of Boolean vectors $\{x_i\}_{i=1}^K \subset \B^N$ is said to be \textit{free} if given any two non-empty subsets $I_1, I_2 \subset \{1, \dots, K\}$ with $I_1 \neq I_2$, we have that 
\[
\bigvee_{I_1} x_i \neq \bigvee_{I_2} x_i.
\]
\end{definition}

We note that freeness is stronger than (linear) independence. Indeed, freeness is a statement about the uniqueness of mixing.  Vectors $\{x_i\}_{i=1}^K \subset \B^N$ are free if and only if each element in $\cone(\{x_i\})$ can be expressed as a linear combination in at most one way.



The next result states the exact criterion for unique $H$ in a Boolean factorization. For that result, we require the following notation. 

\begin{definition}
For $x,y \in  \B^N$, we say that $x$ is \textit{dominated by $y$}, written $x \leq y$, if for each $x_i = 1$ we have that $y_i=1$.  We say $x$ is \textit{strictly dominated by $y$}, written $x < y$ if $x \leq y$ and $x \neq y$.  Given a vector $x \in \B^N$, and a Boolean matrix $W \in \B^{N, R}$, we set
\[
P(x,W) :=\{w_i: w_i \leq x\}.
\]
\end{definition}

If $X = W \land H$, and $x$ is a column of $X$, then $P(x,W)$ are the vectors in $W$ that are dominated by $x$.  That is, the vectors that when added  \textit{may} yield $x$.

\begin{theorem}
\label{thm: unique H if and only if}
Suppose $X \in \B^{N, M}$, $W \in \B^{N,R}$ and $\cone(X) \subset \cone(W)$.  Then $X = W \land H$ has a unique $H$ for the fixed $W$ if and only if the sets $P(x_i, W)$ are free for each column $x_i$. 
\end{theorem}

\begin{proof}
Suppose that each $P(x_i, W)$ are free.  Then any combination of the columns of $W$ from $P(x_i,W)$ yields a unique representation for each $i$.  In particular since $X = W \land H$, for each $i = 1, \dots, M$, 
\[
x_i = \bigvee_{w \in P(x_i, W)} w
\]
is the unique mixing for $x_i$.

Now suppose that $X = W \land H$ has a unique $H$ for $W$.  Suppose to the contrary that there exists some $i$ such that $P(x_i,W)$ is not free. So, there exists $I_1, I_2 \subset P(x_i, W)$ with $I_1 \neq I_2$ non-empty such that 
\[
\bigvee_{w \in I_1} w = \bigvee_{w\in I_2} w.
\]
This yields a contradiction depending on one of two possible cases.

Suppose that $I_1$ is a proper subset of $I_2$. Then since $\vee_{w \in I_1} w = \vee_{w\in I_2} w$,
\[
\begin{array}{rcl}
x_i &=& \bigvee_{w \in P(x_i, W)} w \\
&=& \bigvee_{w \in P(x_i, W)\setminus I_2} w  \bigvee_{w \in I_2} w \\
&=& \bigvee_{w \in P(x_i, W)\setminus I_2} w \bigvee_{w \in I_1} w.
\end{array}
\]
Since $I_1$ is a proper subset of $I_2$, this shows that $x_i$ has two different decompositions, contradicting uniqueness. Similarly, if $I_1$ is not a subset of $I_2$, then
\[
x_i = \bigvee_{w \in P(x_i, W)} w = \bigvee_{w \in P(x_i, W)\setminus I_1} w  \bigvee_{w \in I_2} w
\]
is two different decompositions of $X_i$, also contradicting uniqueness. 
\end{proof}

\begin{corollary}
\label{corollary: unique H}
Suppose $X \in \B^{N, M}$, $W \in \B^{N,R}$ and $\cone(X) \subset \cone(W)$.  If the columns of $W$ are free, then $X = W \land H$ has a unique $H$  for the fixed $W$.
\end{corollary}

\begin{proof}
If the columns of $W$ are free, then $P(X_i,W)$ are free for each $i$.  Indeed, given two non-empty subsets $I_1, I_2 \subset P(x_i,W)$ with $I_1 \neq I_2$, we have that $I_1, I_2 \subset \{w_j\}$.  Hence, $\vee_{I_1} w_j \neq \vee_{I_2} w_j$. 
\end{proof}

As in the nonnegative case,  the converse of Corollary \ref{corollary: unique H} does not hold, as the following example illustrates:

\begin{example}

Let
\[
X = \begin{bmatrix}
0 & 1 & 1 \\
1 & 0 & 1 \\
1 & 1 & 0
\end{bmatrix}
\]
and let $W = X$.  Then the unique $H$ such that $X = W \land H$ is $H = I$, however the columns of $W$ are not free since $x_1 \vee x_2 = x_2 \vee x_3$. 
\end{example}

\subsection{Uniqueness of the decomposition}

It follows from the preceding that in NMF (Boolean factorization), if both Propositions \ref{prop: unique W} and \ref{prop: NMF unique H if and only if} (Propositions \ref{prop: unique W} and Theorem \ref{thm: unique H if and only if}) hold, then both $W$ and $H$ are unique.  Hence,  the corresponding decomposition is unique. However, under an additional rank constraint, the more complicated unique mixing criterion can be dropped. This is what we will explore in this subsection. 
\subsubsection{Uniqueness of NMF}

If $\rk(X) < \rkp(X)$, then NMF is almost never unique.  This occurs because any NMF $X=WH$ solution necessarily utilizes more equations than is needed to span the data $X$. Consequently, almost every point in the interior of $\mbox{cone}(X)$ cannot have a unique representation in terms of the extreme rays \cite{gillis2012sparse}. On the other hand, if $\rkp(X) = \mbox{rk}(X) = R$ then the  subspace spanned by the columns of $X$ has $R$ linearly independent columns.  This forces columns of $W$ to be linearly independent, so that the feature weights must be unique.  In \cite{tam1981geometric}, Tam proves the following:

\begin{theorem}[\cite{tam1981geometric} - Theorem 4.1]
\label{thm: NMF unique}
Let $X \in \R_+^{N,M}$ with $\rk(X) = \rkp(X) = R$.  Then the nonnegative rank factorization $X = WH$, $W \in \R_+^{N,R}$, $H\in \R_+^{R,M}$ is unique if and only if there exists a unique simplicial cone $\C$ such that $\cone(X) \subset \C \subset \mbox{range}(X) \cap \R_+^N$. 
\end{theorem}

\subsubsection{Uniqueness of Boolean factorization}

In NMF, the rank statement $\rk(X) = \rkp(X)$ encodes a linear independence of the columns of $W$. As commented on above, the notion of linear independence does not hold in $\B^N$. Hence, the concept of rank is does not have an exact translation.  However, we have a substitution in $\B^N$ for linear independence, namely freeness. This stronger notion of linear independence naturally leads to a definition of type of column rank for Boolean matrices.
\begin{definition}
Given a matrix $X \in \B^{N, M}$, the \textit{free (Boolean) rank of $X$} denoted $\fr(X)$ is the size of the largest free Boolean subset of $\cone(X) \subset \B^N$. 
\end{definition}

\begin{example}
Consider once again the matrix
$$X = \begin{bmatrix}
1 & 1 & 0 & 0\\
1 & 0 & 1 & 0\\
0 & 1 & 0 & 1\\
0 & 0 & 1 & 1
\end{bmatrix}.$$
From Example \ref{ex: C}, we know that $\rkb(X) = 4$.  However, the largest size of a set of columns one can find that are free is two, so $\fr(X) = 2$.  We once again also note that $X = IX=XI$, so that $X$ fails to have a unique factorization. 
\end{example}

The column rank of a real matrix $X$ is the size of a largest linearly independent subset of $\span(X)$. Similarly, the row rank of $X$ is  the size of a largest linearly independent subset of $\span(X^T)$.  It is well known from elementary linear algebra that the column rank and row rank for real valued matrices are equal to the real rank. While we have defined the free rank of $X$ to correspond to the column rank, we could equally have defined the free rank for the rows of $X$. The next useful result establishes the equality of the free row and column definitions to injectivity:

\begin{proposition}
\label{prop: equiv free rank}

Let $X \in \B^{N, M}$.  Then the following are equivalent:

\begin{enumerate}
    \item $\fr(X) = R$,
    \item The size of the largest subset $S \subset \B^M$ on which $X$ is injective is $R$,
    \item $X$ contains a $R \times R$ permutation submatrix.
\end{enumerate}
\end{proposition}

\begin{proof}
First we recall that $X$ is injective on a set $S \subset \B^M$ if and only if for any pair $z_1, z_2 \in S$, we have 
\[
\bigvee_{i=1}^M z_{i,1} X_i = X \wedge z_1 \neq X \wedge z_2 = \bigvee_{i=1}^M z_{i,2} X_i.
\]
Note the correspondence between $z \in S$ and elements in $\cone(X)$ given by the multiplication $X \wedge z$.  We then see that injectivity on $S$ is equivalent to $\{X \wedge z : z \in S\}$ is free.  Hence, (1) and (2) are equivalent. The equivalence of (2) and (3) follows from Lemma 7.2.6 in \cite{gaubert1992theorie}.



\end{proof}

\begin{corollary}
\label{cor:bool column rank}
Given a matrix $X \in \B^{N, M}$, $\fr(X) \leq \rkb(X)$.
\end{corollary}

\begin{proof}
By Proposition \ref{prop: equiv free rank}, we have that $\fr(X) = R$ if and only if $X$ contains a $R \times R$ permutation matrix.  If $X$ has such a $R\times R$ submatrix, then $\rkb(X) \geq R$ so that $\rkb(X) \geq R = \fr(X)$.
\end{proof}


Proposition \ref{prop: equiv free rank} accomplishes several things.  It connects free rank to the familair concept of injectivity of a matrix on a set.  For Boolean matrices, the largest sets for injectivity correspond to columns that contain the largest permutation matrices, which in turn, bound the free ranks by the Boolean rank. In the case of nonnegativity, we noted that when the matrix failed to be injective (full real rank) one was unlikely to achieve uniqueness of the factorization. The free column rank will play a role similar to real rank  in NMF in Theorem \ref{thm: NMF unique}.  When the free rank is as large as possible, one achieves a unique Boolean matrix factorization:

\begin{theorem}
\label{theorem: unique Boolean NMF}
Let $X \in \B^{N, M}$ and $R = \rkb(X)$.  Then  $\fr(X)  = R$ if and only if every rank $R$ factorization $X=WH$ is permutation equivalent to
\[
W = \begin{pmatrix}
I \\
W' 
\end{pmatrix}
\quad
\quad
\quad
H = \begin{pmatrix}
I & H'
\end{pmatrix},
\]
where $I$ is the $R \times R$ identity matrix and $H' \in \B^{R,M-R}$ and $W' \in \B^{N-R,R}$ are some fixed matrices.  In particular, $X=WH$ has a unique factorization.
\end{theorem}

\begin{proof}
Suppose that every rank $R$ factorization $X=WH$ is permutation equivalent to
\[
W = \begin{pmatrix}
I \\
W' 
\end{pmatrix}
\quad
\quad
\quad
H = \begin{pmatrix}
I & H'
\end{pmatrix}, 
\]
where $I$ is the $R \times R$ identity matrix, $H' \in \B^{R,M-R}$ and $W' \in \B^{N-R,R}$. Multiplying such $W$ and $H$, we see that $X$ is permutation equivalent to
\[
WH = \begin{pmatrix}
I & H'\\
W'& W'H'
\end{pmatrix}.
\]
In particular, $X$ contains an $R \times R$ permutation submatrix.  By Proposition \ref{prop: equiv free rank}, the free rank $\fr(X) = R$. 

Conversely, by Proposition \ref{prop: equiv free rank}, $\fr(X) = R$ if and only if $X$ similar to 
\[
X' = \begin{pmatrix}
I & X_{1,2}\\
X_{2,1} & X_{2,2}
\end{pmatrix}.
\]
where $I$ is the $R \times R$ identity matrix, $X_{1,2} \in \B^{R,M-R}$, $X_{2,1} \in \B^{N-R,R}$, and $X_{2,2} \in \B^{N-R,M-R}$.  Since $\rkb(X)$ is also equal to $R$, we can write $X' = W H$ with
\[
W = \begin{pmatrix}
W_1 \\
W_2 
\end{pmatrix}
\quad
\quad
\quad
H = \begin{pmatrix}
H_1 & H_2
\end{pmatrix},
\]
and $W_1, H_1 \in \B^{R,R}$, $W_2 \in \B^{N-R,R}$, and $H_2 \in \B^{R,M-R}$. Thus we have
\[
\begin{pmatrix}
I & X_{1,2}\\
X_{2,1} & X_{2,2}
\end{pmatrix}
=X'=WH = 
\begin{pmatrix}
W_1 \\
W_2 
\end{pmatrix}
\begin{pmatrix}
H_1 & H_2
\end{pmatrix}
=
\begin{pmatrix}
W_1 H_1 & W_1 H_2\\
W_2 H_1 & W_2 H_2
\end{pmatrix}.
\]
Equating blocks, we see $I = W_1 H_1$. A Boolean matrix is left invertible if and only if it is right invertible if and only if it is invertible \cite{reutenauer1984inversion}.  However, every invertible Boolean matrix must be a permutation matrix \cite{ledley1965inverse}.  Therefore, $W_1 = P$ and $H_1 = P^{-1}$ for some permutation matrix. It then follows that $H_2 = P^{-1} X_{1,2}$ and $W_2 = X_{2,1} P$. That is,
\[
W = \begin{pmatrix}
P \\
X_{2,1} P
\end{pmatrix}
= \begin{pmatrix}
I \\
X_{1,2}
\end{pmatrix}
P
\quad
\mbox{and}
\quad
H = \begin{pmatrix}
P^{-1} & P^{-1} X_{1,2}
\end{pmatrix}
= P^{-1}\begin{pmatrix}
I &  X_{1,2}
\end{pmatrix}.
\]
Thus up to permutation, every factorization of $X=WH$ must have the form
\[
W = \begin{pmatrix}
I \\
W' 
\end{pmatrix}
\quad
\quad
\quad
H = \begin{pmatrix}
I & H'
\end{pmatrix},
\]
\end{proof}

The recent paper \cite{miron2021boolean} appears to be the only other work which contains uniqueness results for Boolean factorization. They characterize uniqueness via a property they called partial uniqueness, namely, given a fixed rank decomposition of $X = \vee_{j=1}^R X_j$, a rank one matrix $X_i$ in the rank decomposition is \textit{partially unique (with respect to the decomposition)} if the only matrix $Y$ that satisfies $X = \vee_{j \neq i} X_j \vee Y$ is $X_i$. The authors of \cite{miron2021boolean} show that the rank decomposition is unique if and only if each rank one matrix in the decomposition is partially unique. Partial uniqueness, and hence their characterization of uniqueness, differs from our work, in the sense that it places a simultaneous constraint on the row and columns of the decomposition $W \land H$. Our consideration for uniqueness are statements about $W$ (uniqueness of cones) and $H$ given a $W$ (freeness of $P(X_i,W)$). In particular, Theorem \ref{thm: unique H if and only if} does not require the user to have acquired the feature weights $H$.  The experimental sections that appear below with regards to stability demonstrate the usefulness of our approach.

\section{Boolean Matrix Factorization with Automatic Model Selection}

An important problem for all factorization methods is the model selection, that is, the estimation of the (usually unknown) number of latent features. Various heuristics to solve this problem have been proposed, including, Akaike’s information criterion (AIC) \cite{akaike1974new}, Bayesian information criterion (BIC) \cite{schwarz1978estimating}, minimum description length (MDL) \cite{rissanen1978modeling}, L-curve method \cite{hansen1992analysis}, and stability method \cite{brunet2004metagenes}.  Automatic Relevance Determination (ARD) method, introduced for neural networks by MakCay \cite{mackay1994automatic}, and applied later for PCA by Bishop \cite{bishop1999bayesian}, and for NMF by Fevotte and Tan \cite{fevotte2009nonnegative} and Morup and Kai \cite{ morup2009tuning} is also popular. Some of the heuristics applied to NMF model selection have been also applied to BMF \cite{miettinen2014mdl4bmf, makhalova2021below}. However, since the nonnegative rank and Boolean rank can be different, using NMF-specific model selection algorithms could potentially produce incorrect estimation for the latent dimension in the Boolean case.

Ideally, one would like a model selection that allows for the identification of unique signatures when present. The following example demonstrates a challenge of this criteria. Namely, that small perturbations of a binary matrix could result in a loss of uniqueness in the Boolean decomposition:

\begin{example}
\label{ex: destroy unique}
Let $X$ be the Boolean matrix given by
\[
X = 
\begin{bmatrix}
1 & 1 & 0 & 0\\
1 & 1 & 0 & 0\\
0 & 0 & 1 & 0\\
0 & 0 & 0 & 1\\
\end{bmatrix}
\]
Then $\rkb(X) = 3$. Moreover, $X$ contains an $3 \times 3$ permutation matrix.  Thus, $\fr(X) = 3$ by Proposition \ref{prop: equiv free rank}.  By Theorem \ref{theorem: unique Boolean NMF}, $X$ has unique Boolean factorization. Now, consider the matrix $Y$ obtained by a single bit flip:
\[
Y = 
\begin{bmatrix}
1 & 1 & 0 & 0\\
1 & 1 & 0 & 0\\
0 & 0 & 1 & \textcolor{red}{1} \\
0 & 0 & 0 & 1\\
\end{bmatrix}
\]
From this minimal perturbation, we no longer have unique decomposition for $Y$.  Indeed, the matrices $W_1$,  $W_2$ given by

\[
W_1 = 
\begin{bmatrix}
1  & 0 & 0\\
1  & 0 & 0\\
0  & 1 & 0 \\
0  & 0 & 1\\
\end{bmatrix}
\quad
W_2 = 
\begin{bmatrix}
1  & 0 & 0\\
1  & 0 & 0\\
0  & 1 & 1 \\
0  & 0 & 1\\
\end{bmatrix}
\]
can both be used to generate the data $Y$.
\end{example}

What Example \ref{ex: destroy unique} shows is that a small perturbation of the original matrix $X$ could result in a loss of uniqueness.  Thus unique factorizations, when they exist, are not stable solutions. In practical applications obtaining stable approximations is often more desirable than unstable exact solutions due to the presence of noise in the data and subsequent interpretations of the decomposition. We therefore desire a model selection algorithm which 1) discovers a stable approximation, such that 2) uniqueness is recovered when it exists.

Here we introduce a heuristic model selection algorithm for estimating the true number, $K$, of original Boolean sources, based on a stability criteria. We refer to this method as BMF$k$. BMF$k$ is an analog for the Boolean semiring of a recent NMF automatic model determination method called  NMF$k$ \cite{alexandrov2013signatures}. It has been shown that when applied to a large number of synthetic datasets with predetermined latent dimensions, NMF$k$ demonstrated a superior performance in comparison to other heuristics \cite{nebgen2021neural}. We aim to seek similar robust results for the Boolean semiring. In this section, we describe the NMF$k$ and BMF$k$ algorithms.

\subsection{Robust Model Selection Algorithm}

\begin{algorithm}

  \caption{Stable Latent Dimension Selection Procedure}
	\KwIn{\\ \(X \in \SR^{M \times N}\): data \\
	\(\mbox{Krange} \subset \mathbb{N}\): a list of explored latent dimensions \(k\) \\ 
	\(S \in \mathbb{N}\): number of random resamplings and decompositions}
	\KwOut{Silhouette scores, $s_k$, relative errors, $e_k$}

	\tcc{Computation part}

	\For{\(k\) in \(\mbox{Krange}\)}{
		\For{\(s\) from \(1\) to \(S\)}{
			$X^{(s)} \sim \mathcal{P}(X)$ \tcp*{Draw Sample from Random Ensemble} \label{alg:k:resamp}
			$W^{(s)}, H^{(s)} \gets \mbox{decompose}(X^{(s)},k)$\tcp*{Decompose Sample of Ensemble} \label{alg:k:decomp}
		}
		
		$C \gets \mbox{CustomCluster}(W^{(1)},\hdots, W^{(S)})$ \tcp*{Cluster Solutions} \label{alg:k:clust}
	    $s_k \gets \mbox{SilhouetteScore}(C)$ \tcp*{Evaluate Cluster Stability} \label{alg:k:sil}
		$e_k \gets \mbox{ReconstructError}(X,C)$ \tcp*{Evaluate Cluster Fit} \label{alg:k:err}
	}
	\label{alg:k}
\end{algorithm}

NMF$k$ and our proposed BMF$k$ are heuristic methods that select latent dimensions of the initial matrix that give stable, or comparatively unique, nonnegative and Boolean approximations respectively. The general procedure for both of these algorithms checks each candidate latent dimension by: decomposing an ensemble of random matrices with mean the initial matrix, clustering the solutions, and evaluating the clusters stability and fit. Algorithm~\ref{alg:k} outlines the unifying procedure that NMF$k$ and BMF$k$ follow to select latent dimensions that provide stable approximations for their respective decompositions.

In Algorithm~\ref{alg:k}, the drawing of matrices from a prescribed random distribution on line~\ref{alg:k:resamp} constructs an ensemble of random matrices, that are slight perturbations of the initial matrix, to mitigate overfitting, to be decomposed with the appropriate algebra on line~\ref{alg:k:decomp}. The details for these two steps vary between NMF$k$ and BMF$k$ and are reported in the following sections. Following their respective random samplings and decompositions, custom clustering is applied in line~\ref{alg:k:clust} to the left factors of the set of decompositions, $\{W^{(1)}, \hdots, W^{(S)}\}$. This custom clustering is an iterative centroid based algorithm similar to $k$-means for NMF$k$ and $k$-medians for BMF$k$. The clustering assigns each column of each $W^{(s)}$ to a centroid, and updates the centroids based on the clusters \cite{nebgen2021neural}. Our custom clustering differs from $k$-means and $k$-medians with the additional constraint that each centroid can only be assigned one column of each $W^{(s)}$. This additional constraint enforces a 1-1 and onto mapping between the columns of each solution $W^{(s)}$ and the current centroids. This constraint simplifies the cluster assignment problem into solving $S$ linear sum assignment problems on distance matrices between each $W^{(s)}$ and the current centroids. After each column assignment is done, the centroids are updated, and the procedure is iterated until converged. 

The evaluation of the resulting clusters is done with two metrics, the silhouette score \cite{rissanen1978modeling} on line~\ref{alg:k:sil} that measures the clusters stability, and a relative error metric on line~\ref{alg:k:err} which measures the quality of fit. Silhouette scores range between -1, the poorest quality clustering, to 1, the highest quality clustering. We aggregate the silhouettes into a single score by taking the minimum of the average of the silhouettes for each cluster. A minimum silhouette score close to 1 indicates that the solutions of each $X^{(s)}$ in the ensemble form tight clusters. The relative error measure is evaluated from the medoids of the clusters and an $H$ computed through nonnegative or Boolean regression. Since this process is done for each latent dimension of interest, L-statistics \cite{vangara2021finding} can then be used to determine the optimal stable dimension, the elbows of the reconstruction errors and the silhouette scores. Further details of resampling, decomposing, and clustering for the individual implementations of NMF$k$ and BMF$k$ are described in the following sections.

\subsection{NMF$k$}

NMF$k$ discovers the latent dimension that gives a robust nonnegative approximation for a given nonnegative matrix $X$ by following the procedure given in Algorithm~\ref{alg:k}. The NMF$k$ algorithm begins by drawing samples from a random distribution to form a random ensemble. While there is an abundance of suitable distributions to draw from for nonnegative data, our reported results use the elementwise resampling defined by $\mathcal{P}(x) = x \cdot \mathcal{U}(1-\epsilon, 1+\epsilon)$ where $\mathcal{U}(a,b)$ is a uniform distribution on the interval $[a,b)$. Alternative noise models can be selected appropriately for various matrices and applications.

For each element of the random ensemble, the sample is decomposed according to a prescribed noise model and decomposition algorithm on line~\ref{alg:k:decomp}. There are numerous NMF objective functions corresponding to different noise models, e.g. Kullback-Liebler divergence, Frobenius norm, Itakura-Saito divergence, as well as numerous algorithms for each objective function e.g. Multiplicative Update, Alternating Direction Method of Multipliers, Block Principal Pivoting \cite{cichocki2009nonnegative}. NMFk decomposes each matrix in the ensemble, $X^{(s)}$, using a multiplicative update algorithm  \cite{lee1999learning} to solve the Frobenius norm NMF problem:
\begin{equation}
\begin{aligned}
& \argmin_{W,H}
& & \frac{1}{2}||X^{(s)}-W H||^2_F \\
& \text{subject to} & & W \in \R_+^{N \times k},\\
& & & H \in \R_+^{k \times M}.
\end{aligned}
\label{eqn: NMF problem}
\end{equation}

Multiplicative updates is an iterative algorithm that alternates optimizing over $W$ and $H$ with the update rules 
\begin{equation}
\begin{aligned}
    W &\gets W \frac{XH^\top}{WHH^\top}\\
    H &\gets H \frac{W^\top X}{W^\top W H}
\end{aligned}
\end{equation}
always using the updated variables in the subsequent computation. This procedure is iterated until a fixed number of steps, or some convergence criteria is met. The multiplicative update algorithm preserves nonnegativity given nonnegative initializations, which were taken to be uniform random matrices. 

In the clustering step of NMF$k$ we employed a cosine similarity  in the linear sum assignment problem to determine cluster membership. After the cluster assignments, the cluster centroids were updated as the medians of the assigned members to remove outlying solutions. Similarly, a cosine distance was employed in the silhouette statistic computation. In practice, the noise resampling technique, the NMF objective function and algorithm, and the distance metrics can all be selected for the application and the NMF$k$ procedure successfully aids in the selection of a latent dimension that provides stable solutions.

\subsection{BMF$k$}


The BMF$k$ algorithm follows the series of steps in Algorithm~\ref{alg:k}, adjusted for Boolean factorizations. Once again, an ensemble of BMF solutions $(W_s, H_s)_{s=1}^S$ is obtained from an ensemble of perturbed matrices $\{X^{(1)}, \hdots, X^{(S)}\}$ for each permissible latent dimension $k$. 
For BMF$k$, the ensemble of random matrices is generated by randomly flipping elements of the original matrix $X$ according a Bernoulli distribution. This corresponds to resampling the matrix $X$ using the element-wise distribution $\mathcal{P}(x) = x +_f \mathcal{B}(p)$, where the $+_f$ operator is the bit flipping operator, and $\mathcal{B}(p)$ is the Bernoulli distribution with parameter $p$. 
The sampled matrices are then decomposed into a Boolean factorization using one of many BMF algorithms. In this work our reported results use the BANMF algorithm \cite{truong2021boolean} with random initial conditions.For recent developments in the area of BMF, please see the review \cite{miettinen2020recent}. The BANMF algorithm solves the Boolean matrix factorization problem for each $X^{(s)}$ through the auxiliary optimization problem:
\begin{equation}
\begin{aligned}
& \argmin_{Y,W,H}
& & ||Y-W H||_F \\
& \text{subject to}
& & 1 \leq Y_{ij} \leq k, \text{ if } X^{(s)}_{ij} = 1 \\
& & & Y_{ij} = 0, \text{ if } X^{(s)}_{ij} = 0 \\
& & & W \in \R_+^{N \times k},\\
& & & H \in \R_+^{k \times M}.
\end{aligned}
\label{eqn: BANMF problem}
\end{equation}

The factors \(W\) and \(H\) are updated via an alternating optimzation procedure similar to the multiplicative updates, while the \textit{auxiliary matrix \(Y\)} is updated by being projected to the constraint set according to the update rules,
\begin{equation}
\begin{aligned}
    W &\gets W \frac{YH^\top}{WHH^\top}\\
    H &\gets H \frac{W^\top Y}{W^\top W H}\\
    Y &\gets X^{(s)}\star \operatorname{maximum}(\operatorname{minimum}(WH,k),1),
\end{aligned}
\end{equation}
where $\star$, $\operatorname{maximum}$, and $\operatorname{minimum}$ are the elementwise multiplication, maximum operator, and minimum operator respectively.
The solution of this optimization problem will then be thresholded according to the following minimization problem,
\begin{equation}
\begin{aligned}
    \min_{w \in \R, h \in \R} ||X^{(s)} - (W \geq w) \otimes_B (H \geq h) ||_F
\end{aligned}
\end{equation}
to arrive at a BMF solution. The grid-search thresholding algorithm is shown in \cite{truong2021boolean}.

In the clustering step, since the factors are Boolean, the custom clustering algorithm uses hamming distance to determine cluster assignment for each of the $K$ features of each $W^{(s)}$. The centroids are updated by taking the elementwise medians of the vectors assigned to each cluster as this minimizes the hamming distance objective function, $\min_{c_k} \sum_{w \in \mathcal{W}_k} \operatorname{hamming}(c_k, w)$, to find a centroid $c_k$ for a cluster $\mathcal{W}_k = \{ W^{(1)}_k, W^{(2)}_k, \hdots, W^{(S)}_k\}$. Silhouette scores with hamming distance are used to determine how well each latent dimension $K$ was clustered.



		
	

\subsection{Experimental evaluation} 

We demonstrate the effectiveness of BMF$k$ at discovering the stable latent dimension in several synthetically generated datasets.  Our first example compares the extracted signals from BMF$k$ to NMF$k$ on a matrix with different nonnegative and Boolean ranks.  For our second example, we generate a large collection of Boolean matrices with unique BMF, and show that not only does BMF$k$ identify the correct latent dimension, but also the unique factors. We then inject noise into a random collection of Boolean matrices with unique BMF and show that BMF$k$ again discovers the correct latent dimension, and the discovered features are highly correlated with the ``true'' hidden features. 

To evaluate an extracted feature, $a_1$, relative to a corresponding generative feature, $a_2$, we rely on the cosine similarity metric, $$\dfrac{a_1^\top \cdot a_2}{||a_1||||a_2|}\;.$$ The cosine similarity applied to all pairs of recovered and generative features constructs a cosine similarity matrix. 

To evaluate a decomposition, we aggregate the cosine similarities of both the left and right factors into a single \textit{score metric} to measure how well the ground-truth decomposition is recovered \cite{battaglino2018practical}. For a pair of rank-one matrices \(X_1 = a_1b_1^\top\) and \(X_2 = a_2b_2^\top\), the score metric is defined as:
\[\text{score}(X_1,X_2) = \dfrac{a_1^Ta_2}{||a_1||||a_2||}  \dfrac{b_1^Tb_2}{||b_1||||b_2||}\;. \]
For higher rank matrices and decompositions, the average of the scores of all rank one factors is taken, after the decompositions are permuted to maximize the scores. 
For reference, we also report results of other Boolean model selection strategies: (1) Minimum Description Length (MDL) \cite{rissanen1978modeling} and (2) Covarge Quality (CQ) \cite{trnecka2021model}. While BMF$k$ focuses on the solution stability, MDL minimizes the description of the data using the BMF model and CQ measures the change of the angle in the coverage (or reconstruction) curve.
\subsubsection{Dataset with $\rkb(X) \neq \rkp(X)$}

\begin{figure}[t]
	\centering
	\includegraphics[width = 0.8 \textwidth]{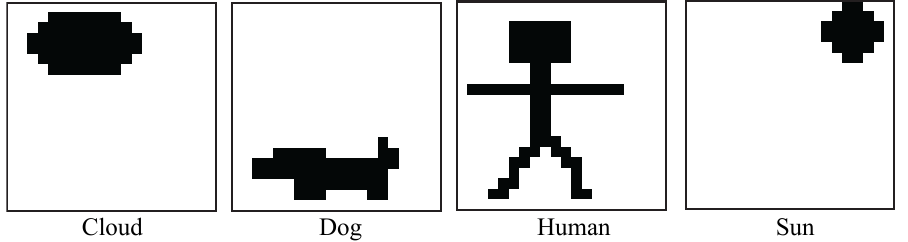}
	\caption{Four binary images are used as Boolean latent features to generate the synthetic image data.}
	\label{fig:imagedata}
\end{figure}

This example demonstrates that the underlying semiring can affect the latent dimension in practice. A synthetic image dataset is generated whose latent features are four binary \(20 \times 20\) images, as shown in Figure \ref{fig:imagedata}. All possible (nonempty) Boolean combination images are generated, producing a dataset $X$ of size \(400 \times 15\) with Boolean rank \(\rkb(X) =4\). The overlap between the images `human' and `dog', and between `cloud' and `human' forces \(\rkp(X)=6\). Each image has an entry with a value of $1$ where each other has a value of $0$.  It follows that $X$ contains a $4 \times 4$ permutation matrix, so that by Theorem \ref{theorem: unique Boolean NMF}, $X$ has a unique BMF. 

\begin{figure}[ht]
	\centering
	\includegraphics[width = 0.8 \textwidth]{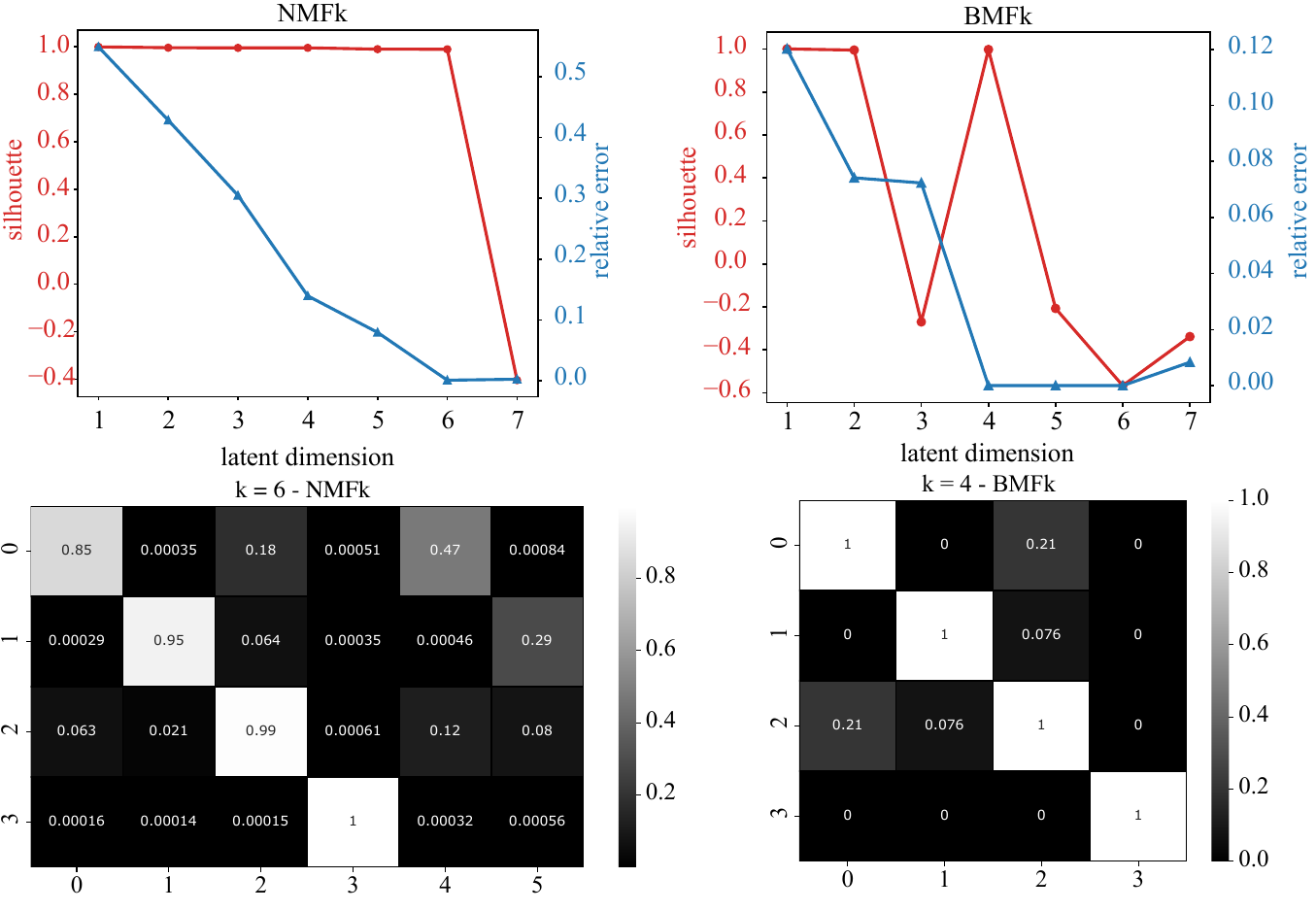}
	\caption{\textit{Top -} NMF$k$ and BMF$k$ silhouette scores and relative errors for candidate latent dimensions. \textit{Bottom -} Cosine similarity matrices comparing the extracted features to the generative Boolean features. }
	\label{fig:imageresult}
\end{figure}

Figure~\ref{fig:imageresult} (top) depicts the resulting silhouette and relative error scores of NMF$k$ and BMF$k$ for each candidate latent dimension. NMF$k$ correctly identifies the nonnegative rank at \(K=6\) and BMF$k$ correctly identifies the Boolean rank of the data of \(K=4\). Figure~\ref{fig:imageresult} (bottom) additionally depicts the score metric between the recovered factors of each method, and the Boolean images used to generate the matrix. Clearly, BMF$k$ perfectly recovered the factors with the scores of 1, while several of NMF$k$'s scores are lower. 

\begin{figure}
	\centering
	\includegraphics[width=0.8 \textwidth]{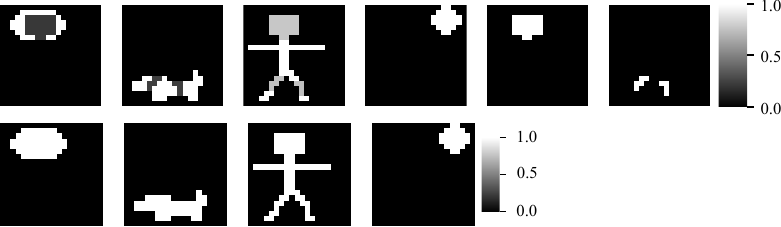}
	\caption{\textit{Top -} NMF$k$ extracted factors. \textit{Bottom -} BMF$k$ extracted factors are perfectly correlated with the true unique factors.}
	\label{fig:imagefactor}
\end{figure}

Figure \ref{fig:imagefactor} shows the extracted images from both methods. Note that the stable features identified from NMF$k$ include the overlapping regions. This is a more robust solution than say, completely separated features where small perturbations could disrupt such signals for reconstruction. The BMF$k$ technique has correctly identified the four unique signals. 

Table~\ref{table:method compare} compares the selected dimension of several Boolean latent dimension selection models on the image data. Clearly, only BMF$k$ and BANMF-MDL correctly identified the latent dimension, while BANMF-CQ selected a lower latent dimension corresponding to a higher relative error plateau seen at $k=2$ and $k=3$.

\begin{table*}[ht]
  \centering
	\begin{tabular}{|c|c|c|c|c|c|}
		\hline
		Datasets & Boolean rank & BMF$k$ & BANMF-MDL & BANMF-CQ \\
		\hline
		image & 4 & 4 & 4 & 2 \\
		\hline
	\end{tabular}
\caption{Comparisons of latent Boolean dimension selection algorithms applied to the image dataset.}
	\label{table:method compare}
\end{table*}

\subsubsection{BMF$k$ stable rank finds unique features}

The next experiments demonstrate that for matrices generated with unique Boolean features, BMF$k$ identifies the latent dimension corresponding to these unique decompositions, and extracts these unique patterns.  We first show that for a class of matrices $X$ with unique BMF, BMF$k$ exactly finds the unique features.  We then add noise to these matrices and show that the BMF$k$ algorithm successfully recovers the correct latent dimension and features that are highly correlated with the original pre-corrupted unique features. 

\begin{figure}[ht]
	\centering
	\includegraphics[width = 0.6 \textwidth]{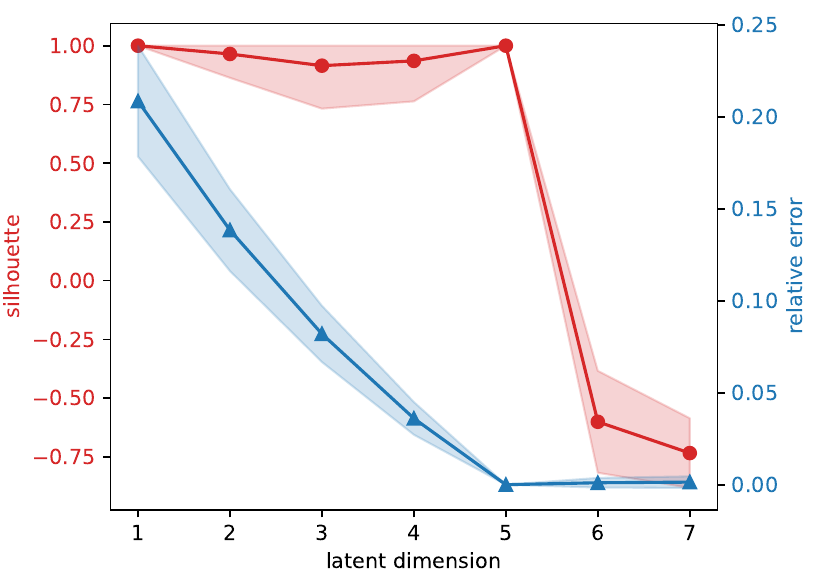}
	\caption{The mean values with one standard deviation band of silhouette scores and Boolean relative error of 1600 random matrices that have a unique factorization at \(\rkb=5\). BMF$k$ correctly identify the Boolean rank, and extracted the unique features for all matrices.}
	\label{fig:random example}
\end{figure}

The dataset used for this experiment is generated  based on the uniqueness result from Theorem \ref{theorem: unique Boolean NMF}. A single data matrix \(X\) is generated as following: First, a \(W' \in \mathcal{B}^{15,5}\) and \(H' \in \mathcal{B}^{5,15}\) with a density of 0.3 are randomly generated. Then \(X = \begin{pmatrix}	I \\W'\end{pmatrix}\)\(\begin{pmatrix}I & H'\end{pmatrix}\). We apply BMF$k$ on 1600 randomly generated dataset and collect the silhouette scores, relative Boolean error, and the score between the extracted features and the true factors. Figure \ref{fig:random example} shows the average with one standard deviation band of silhouette scores and relative Boolean error of 1600 random matrices. Here we see that BMF$k$ successfully identified the correct Boolean latent dimension at \(k=5\) for all matrices. Moreover, the factor scores are ones for all matrices indicating that BMF$k$ can perfectly recover the original factors.

\begin{figure}[ht]
	\centering
	\includegraphics[width = 1.0 \textwidth]{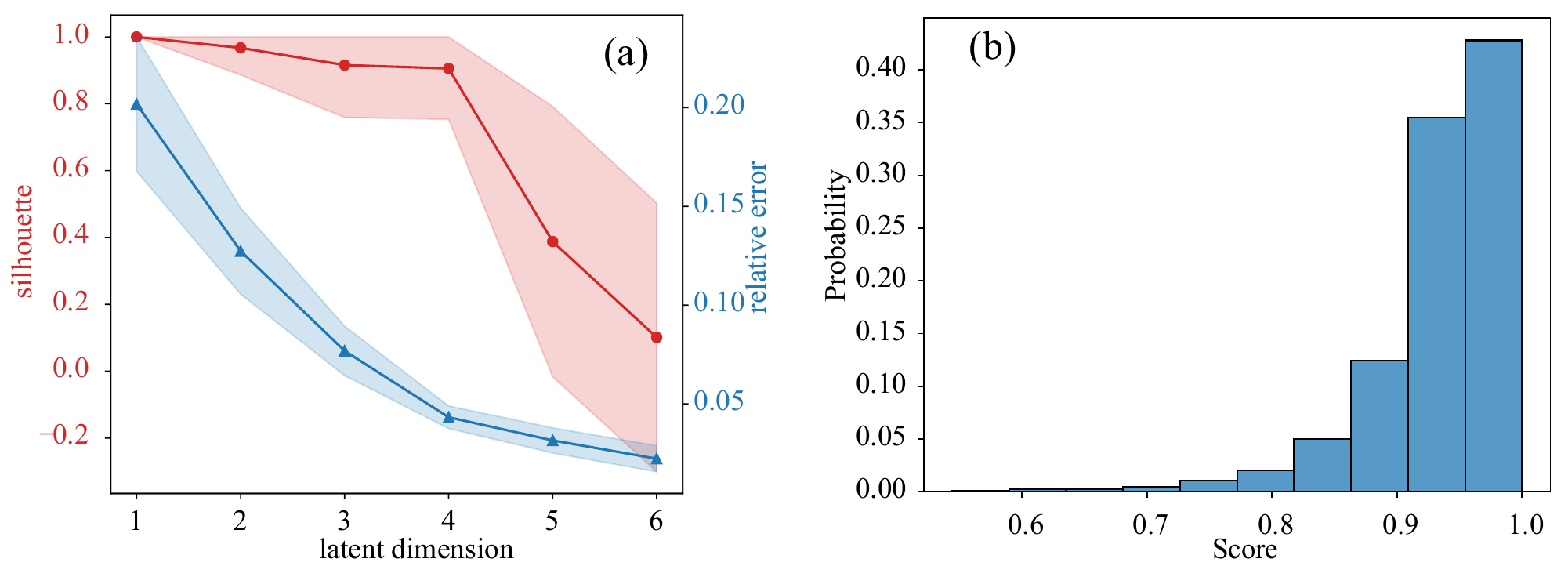}
	\caption{The mean values with one standard deviation band of silhouette scores and Boolean relative error of 1600 random datasets \textbf{with 5\% flipping noise} that have a unique factorization at \(\rkb=4\).}
	\label{fig:random noise 5}
\end{figure}

Next, we generate a dataset using \(W' \in \mathcal{B}^{15,4}\), \(H' \in \mathcal{B}^{4,15}\), density of 0.3 and \(X = \begin{pmatrix}	I \\W'\end{pmatrix}\)\(\begin{pmatrix}I & H'\end{pmatrix}\).  We then take this collection of matrices with unique BMF and corrupt the data by applying modest Bernoulli noise.  Here, $5\%$ of the entries of each matrix $X$ were flipped. Once again, we have plotted the average and standard deviation of the silhouette scores with relative error in Figure \ref{fig:random noise 5} (a). In Figure \ref{fig:random noise 5} (b), we have plotted the distribution of the metric scores from all datasets. In this case, where the input data is noisy, BMF$k$ can extract the ground truth latent features quite well.

\FloatBarrier
\section*{Acknowledgements}

This work was supported by the LDRD program of Los Alamos
National Laboratory under project number 20190020DR and the Center for Nonlinear Studies. Los Alamos National Laboratory is operated by Triad National Security, LLC, for the National Nuclear Security Administration
of U.S. Department of Energy (Contract No. 89233218CNA000001).

\section{Appendix}

\subsection{Proofs for rank results}
Here we provide detailed proofs of Proposition \ref{prop: NMF to Boolean} and Theorem \ref{thm: z2 rank}:

\begin{customthm2}{1}
Let  $X \in \BB^{N,M}$.  If $X = WH$ is a nonnegative factorization of $X$, then $\hat{W} \land \hat{H}$ is a Boolean factorization of $X$. 
\end{customthm2}

\begin{proof}
Note that from the nonnegative factorization,
\[
 \sum_{k=1}^K W_{i,k} H_{k,j} = X_{i,j} \in \{0,1\}.
\]
Because each entry in the summand is nonnegative, the only way the above sum results in zero is if each element is itself zero. Similarly, the only way the above sum results in one is if at least one of the products $W_{i,k} H_{k,j}$ is non-zero.\\ 
Now consider the Boolean product $Y = \hat{W} \land \hat{H}$.   We will show that $Y = X$. By definition,
\[
Y_{i,j} = \bigvee_{k=1}^K \hat{W}_{i,k} \wedge \hat{H}_{k,j}.
\]
If $X_{i,j} = 0$, then by above each $ W_{i,k} H_{k,j} = 0$, and hence each $\hat{W}_{i,k} \wedge \hat{H}_{k,j}=0$.  In this case, $Y_{i,j} = 0$.  Similarly, if $X_{i,j} = 1$, then because there exists some product $ W_{i,k} H_{k,j} \neq 0$, we have that $ \hat{W}_{i,k} \wedge \hat{H}_{k,j} = 1$.  In this case, $Y_{i,j} = 1$.  Therefore, we have $Y = X$.
\end{proof}

Next, we show that if $X \in \BB^{N,M}$, then $\rkz(X) \leq \rk(X)$. Recall that if $X$ is an $R \times R$ matrix over a field $\mathbb{F}$, one can define its determinant. Concretely, if $X = (X_{i,j})_{i,j=1}^R$, then 
\[
\mbox{det}(X):= \sum_{\sigma \in S_R} \mbox{sgn}(\sigma) \prod_{i=1}^R X_{\sigma(i),i},
\]
where $\sigma$ is a permutation of $\{1, \dots, R\}$,  $S_R$ is the set of all permutations, and $\mbox{sgn}(\sigma) \in \{-1, 1\}$ is the sign of the permutation $\sigma$. 

If $X$ is an $R \times R$ matrix over a field $\mathbb{F}$, then $\mbox{det}(X) \neq 0$ if and only if $X$ is invertible \cite{lang1993algebra}. For a  $N \times M$ matrix $X$ over a field $\mathbb{F}$, the row rank, column rank, and factor rank are all equal \cite{roman2005advanced}. Therefore  $X$ has rank $R$ if and only if $X$ contains a $R \times R$ submatrix $X'$ with full rank, and every larger square submatrix is not full rank. Thus $X \in \mathbb{F}^{N,M}$ has rank $R$ if and only if $X$ contains a $R \times R$ submatrix $X'$ with $\mbox{det}(X') \neq 0$ and every larger square submatrix has determinant zero.

\begin{customthm}{2}
Let $X \in \BB^{N,M}$. Then $\rkz(X) \leq \rk(X)$.
\end{customthm}

\begin{proof}
Let $\rk(X) = R$.  By the comments proceeding the theorem, if $K > R$, then every $K \times K$ submatrix $X'$ of $X$ will have real determinant $0$. That is, 
\[
\mbox{det}_\R(X')= \sum_{\sigma \in S_K} \mbox{sgn}(\sigma) \prod_{i=1}^K X_{\sigma(i),i} = 0
\]
For each $\sigma$, $\prod_{i=1}^K X_{\sigma(i),i}$ is either $0$ or $1$. Hence, the above summation is a some quantity of $1$'s perfectly balanced by an equal quantity of $-1$'s. Therefore, the summation consists of an even number of non-zero elements.

Now, consider the $\mathbb{Z}_2$ determinant of $X'$.  Since $-1 = 1$ in $\mathbb{Z}_2$, $\mbox{sgn}(\sigma)  = 1$ for each $\sigma$. Therefore, we have that 
\[
\mbox{det}_{\mathbb{Z}_2}(X'):= \bigoplus_{\sigma \in S_K} \bigotimes_{i=1}^K X_{\sigma(i),i},
\]
where now the sum and product are happening over $\mathbb{Z}_2$. Note the product $\bigotimes_{i=1}^K X_{\sigma(i),i}$ happening within $\mathbb{Z}_2$ is the same as within $\R$,  and so remains unchanged. From above, we know there is an even number of non-zero elements in this sum.  Hence $\mbox{det}_{\mathbb{Z}_2}(X') = 0$. By the comments preceding the theorem, $\rkz(X) \leq R = \rk(X)$. 

\end{proof}

\subsection{Proof that minimal generating set for Boolean cone unique}

We will show that a minimal generating set for a cone $\C$ over the Boolean semiring $\B^N$ is necessarily unique. To prove this, we require the following order notation:   
\begin{definition}
\label{def: order}
 Let 
\[
[0,y) := \{x \in \C: x < y \}
\]
be the set of vectors strictly dominated by $y$. We say that $y$ is \textit{minimal in $\C$} if 
\[
\bigvee_{x \in [0,y)} x  < y.
\]
Finally, we define 
\[
\min(\C) := \{y \in \C: y \mbox{ is minimal in } \C, y \neq 0\}.
\]
\end{definition}
For example, consider the cone
\[
\C = 
\left\{
\begin{pmatrix}
0 \\
0 \\
0 \\
\end{pmatrix},
\begin{pmatrix}
1 \\
0 \\
0 \\
\end{pmatrix},
\begin{pmatrix}
0 \\
1 \\
0 \\
\end{pmatrix},
\begin{pmatrix}
1\\
0 \\
1 \\
\end{pmatrix},
\begin{pmatrix}
1 \\
1 \\
0 \\
\end{pmatrix},
\begin{pmatrix}
1 \\
1 \\
1 \\
\end{pmatrix}
\right\}
\]
Then $\min(\C) = \{(1,0,0), (0,1,0), (1,0,1)\}$.  We will show that $\min(\C)$ is the unique minimal generating set, and therefore the set of extreme rays are unique.

\begin{theorem}
Let $\C \subset \B^N$ be a cone.  Then the set of extreme rays are unique, and are equal to $\min(\C)$. 
\end{theorem}

\begin{proof}
We will show that $\min(\C)$ is a subset of any generating set. We will then show that $\min(\C)$ also generates $\C$.  It will then follow that the unique minimal generating set is $\min(\C)$.

First, we note that if $G \subset \C$ is a generating set, then $\min(\C) \subset \G$. By definition, if $y \in \min(\C)$ then $y \in \C$ and nothing else from $\C$ can be added to get $y$. Hence, it must be in any generating set $G$.

Next, we show that $\min(\C)$ spans $\C$.  Let $\C' = \mbox{span}_\vee (\min(\C))$. Suppose to the contrary that there exists some $y \in \C$ but $y \notin \C'$.  Then certainly, $y \notin \min(\C)$ so that
\[
\bigvee_{x \in [0,y)} x  = y.
\]
But consider now $[0,y)$. Since $y \notin \C'$, there must exist a $y_1 \in [0,y)$ such that $y_1 \notin \C'$ and
\[
\bigvee_{x \in [0,y) \setminus \{y_1\}} x  < y.
\]
One can then apply this process again on the new vector $y_1$. Since any vector in $\B^N$ has at most $N$ non-zero entries, this process must terminate after a finite number steps $k$. This means that $y_k \notin \C'$, but yet any sub-vector that is required to build $y_k$ must belong to $\C'$.  This is clearly a contradiction, implying that $y$ must have belonged to $\C'$.
\end{proof}

\bibliographystyle{plain}
\bibliography{references.bib}

\begin{thebibliography}{10}

\bibitem{akaike1974new}
Hirotugu Akaike.
\newblock A new look at the statistical model identification.
\newblock {\em IEEE transactions on automatic control}, 19(6):716--723, 1974.

\bibitem{alexandrov2013signatures}
Ludmil~B Alexandrov, Serena Nik-Zainal, David~C Wedge, Samuel~AJR Aparicio, Sam
  Behjati, Andrew~V Biankin, Graham~R Bignell, Niccolo Bolli, Ake Borg,
  Anne-Lise B{\o}rresen-Dale, et~al.
\newblock Signatures of mutational processes in human cancer.
\newblock {\em Nature}, 500(7463):415--421, 2013.

\bibitem{battaglino2018practical}
Casey Battaglino, Grey Ballard, and Tamara~G Kolda.
\newblock A practical randomized cp tensor decomposition.
\newblock {\em SIAM Journal on Matrix Analysis and Applications},
  39(2):876--901, 2018.

\bibitem{beasley2012isolation}
LeRoy~B Beasley.
\newblock Isolation number versus boolean rank.
\newblock {\em Linear algebra and its applications}, 436(9):3469--3474, 2012.

\bibitem{bishop1999bayesian}
Christopher~M Bishop.
\newblock Bayesian pca.
\newblock {\em Advances in neural information processing systems}, pages
  382--388, 1999.

\bibitem{brown1968invertibly}
Thomas~Andrew Brown, Mario~Leon Juncosa, and VL~Klee.
\newblock Invertibly positive linear operators on spaces of continuous
  functions.
\newblock Technical report, RAND CORP SANTA MONICA CALIF, 1968.

\bibitem{brunet2004metagenes}
Jean-Philippe Brunet, Pablo Tamayo, Todd~R Golub, and Jill~P Mesirov.
\newblock Metagenes and molecular pattern discovery using matrix factorization.
\newblock {\em Proceedings of the national academy of sciences},
  101(12):4164--4169, 2004.

\bibitem{cichocki2009nonnegative}
Andrzej Cichocki, Rafal Zdunek, Anh~Huy Phan, and Shun-ichi Amari.
\newblock {\em Nonnegative matrix and tensor factorizations: applications to
  exploratory multi-way data analysis and blind source separation}.
\newblock John Wiley \& Sons, 2009.

\bibitem{cohen1993nonnegative}
Joel~E Cohen and Uriel~G Rothblum.
\newblock Nonnegative ranks, decompositions, and factorizations of nonnegative
  matrices.
\newblock {\em Linear Algebra and its Applications}, 190:149--168, 1993.

\bibitem{doherty1999biclique}
Faun~CC Doherty, J~Richard Lundgren, and Daluss~J Siewert.
\newblock Biclique covers and partitions of bipartite graphs and digraphs and
  related matrix ranks of $\{$0, 1$\}$-matrices.
\newblock {\em Congressus Numerantium}, pages 73--96, 1999.

\bibitem{fevotte2009nonnegative}
C{\'e}dric F{\'e}votte and A~Taylan Cemgil.
\newblock Nonnegative matrix factorizations as probabilistic inference in
  composite models.
\newblock In {\em 2009 17th European Signal Processing Conference}, pages
  1913--1917. IEEE, 2009.

\bibitem{gaubert1992theorie}
St{\'e}phane Gaubert.
\newblock {\em Th{\'e}orie des syst{\`e}mes lin{\'e}aires dans les
  dio{\"\i}des}.
\newblock PhD thesis, Paris, ENMP, 1992.

\bibitem{geerts2004tiling}
Floris Geerts, Bart Goethals, and Taneli Mielik{\"a}inen.
\newblock Tiling databases.
\newblock In {\em International conference on discovery science}, pages
  278--289. Springer, 2004.

\bibitem{gillis2012sparse}
Nicolas Gillis.
\newblock Sparse and unique nonnegative matrix factorization through data
  preprocessing.
\newblock {\em The Journal of Machine Learning Research}, 13(1):3349--3386,
  2012.

\bibitem{gondran2008graphs}
Michel Gondran and Michel Minoux.
\newblock {\em Graphs, dioids and semirings: new models and algorithms},
  volume~41.
\newblock Springer Science \& Business Media, 2008.

\bibitem{gutch2012ica}
Harold~W Gutch, Peter Gruber, Arie Yeredor, and Fabian~J Theis.
\newblock Ica over finite fields—separability and algorithms.
\newblock {\em Signal Processing}, 92(8):1796--1808, 2012.

\bibitem{hansen1992analysis}
Per~Christian Hansen.
\newblock Analysis of discrete ill-posed problems by means of the l-curve.
\newblock {\em SIAM review}, 34(4):561--580, 1992.

\bibitem{lang1993algebra}
Serge Lang.
\newblock Algebra. 3rd.
\newblock {\em Edition Addison--Wesley}, 1993.

\bibitem{laurberg2008theorems}
Hans Laurberg, Mads~Gr{\ae}sb{\o}ll Christensen, Mark~D Plumbley, Lars~Kai
  Hansen, and S{\o}ren~Holdt Jensen.
\newblock Theorems on positive data: On the uniqueness of nmf.
\newblock {\em Computational intelligence and neuroscience}, 2008, 2008.

\bibitem{ledley1965inverse}
Robert~S Ledley.
\newblock The inverse of a boolean matrix.
\newblock Technical report, NATIONAL BIOMEDICAL RESEARCH FOUNDATION WASHINGTON
  DC, 1965.

\bibitem{lee1999learning}
Daniel~D Lee and H~Sebastian Seung.
\newblock Learning the parts of objects by non-negative matrix factorization.
\newblock {\em Nature}, 401(6755):788--791, 1999.

\bibitem{li2005general}
Tao Li.
\newblock A general model for clustering binary data.
\newblock In {\em Proceedings of the eleventh ACM SIGKDD international
  conference on Knowledge discovery in data mining}, pages 188--197, 2005.

\bibitem{mackay1994automatic}
David~JC MacKay and Radford~M Neal.
\newblock Automatic relevance determination for neural networks.
\newblock In {\em Technical Report in preparation}. Cambridge University, 1994.

\bibitem{makhalova2021below}
Tatiana Makhalova and Martin Trnecka.
\newblock From-below boolean matrix factorization algorithm based on mdl.
\newblock {\em Advances in Data Analysis and Classification}, 15(1):37--56,
  2021.

\bibitem{miettinen2020recent}
Pauli Miettinen and Stefan Neumann.
\newblock Recent developments in boolean matrix factorization.
\newblock {\em arXiv preprint arXiv:2012.03127}, 2020.

\bibitem{miettinen2014mdl4bmf}
Pauli Miettinen and Jilles Vreeken.
\newblock Mdl4bmf: Minimum description length for boolean matrix factorization.
\newblock {\em ACM transactions on knowledge discovery from data (TKDD)},
  8(4):1--31, 2014.

\bibitem{miron2021boolean}
Sebastian Miron, Mamadou Diop, Anthony Larue, Eddy Robin, and David Brie.
\newblock Boolean decomposition of binary matrices using a post-nonlinear
  mixture approach.
\newblock {\em Signal Processing}, 178:107809, 2021.

\bibitem{morup2009tuning}
Morten M{\o}rup and Lars~Kai Hansen.
\newblock Tuning pruning in sparse non-negative matrix factorization.
\newblock In {\em 2009 17th European Signal Processing Conference}, pages
  1923--1927. IEEE, 2009.

\bibitem{nebgen2021neural}
Benjamin~T Nebgen, Raviteja Vangara, Miguel~A Hombrados-Herrera, Svetlana
  Kuksova, and Boian~S Alexandrov.
\newblock A neural network for determination of latent dimensionality in
  non-negative matrix factorization.
\newblock {\em Machine Learning: Science and Technology}, 2(2):025012, 2021.

\bibitem{orlin1977contentment}
James Orlin et~al.
\newblock Contentment in graph theory: covering graphs with cliques.
\newblock In {\em Indagationes Mathematicae (Proceedings)}, volume~80, pages
  406--424. North-Holland, 1977.

\bibitem{reutenauer1984inversion}
Christophe Reutenauer and Howard Straubing.
\newblock Inversion of matrices over a commutative semiring.
\newblock {\em Journal of Algebra}, 88(2):350--360, 1984.

\bibitem{rissanen1978modeling}
Jorma Rissanen.
\newblock Modeling by shortest data description.
\newblock {\em Automatica}, 14(5):465--471, 1978.

\bibitem{roman2005advanced}
Steven Roman, S~Axler, and FW~Gehring.
\newblock {\em Advanced linear algebra}, volume~3.
\newblock Springer, 2005.

\bibitem{schwarz1978estimating}
Gideon Schwarz.
\newblock Estimating the dimension of a model.
\newblock {\em The annals of statistics}, pages 461--464, 1978.

\bibitem{stewart1993early}
Gilbert~W Stewart.
\newblock On the early history of the singular value decomposition.
\newblock {\em SIAM review}, 35(4):551--566, 1993.

\bibitem{tam1981geometric}
Bit-Shun Tam.
\newblock A geometric treatment of generalized inverses and semigroups of
  nonnegative matrices.
\newblock {\em Linear Algebra and its Applications}, 41:225--272, 1981.

\bibitem{trnecka2021model}
Martin Trnecka and Marketa Trneckova.
\newblock Model order selection for approximate boolean matrix factorization
  problem.
\newblock {\em Knowledge-Based Systems}, page 107184, 2021.

\bibitem{truong2021boolean}
Duc~P. Truong, Erik Skau, Derek Desantis, and Boian Alexandrov.
\newblock Boolean matrix factorization via nonnegative auxiliary optimization.
\newblock {\em IEEE Access}, 9:117169--117177, 2021.

\bibitem{vangara2021finding}
Raviteja Vangara, Manish Bhattarai, Erik Skau, Gopinath Chennupati, Hristo
  Djidjev, Thomas Tierney, James~P Smith, Valentin~G Stanev, and Boian~S
  Alexandrov.
\newblock Finding the number of latent topics with semantic non-negative matrix
  factorization.
\newblock {\em IEEE Access}, 2021.

\bibitem{yeredor2011independent}
Arie Yeredor.
\newblock Independent component analysis over galois fields of prime order.
\newblock {\em IEEE Transactions on Information Theory}, 57(8):5342--5359,
  2011.

\bibitem{zhang2007binary}
Zhongyuan Zhang, Tao Li, Chris Ding, and Xiangsun Zhang.
\newblock Binary matrix factorization with applications.
\newblock In {\em Seventh IEEE International Conference on Data Mining (ICDM
  2007)}, pages 391--400. IEEE, 2007.

\end{thebibliography}

\end{document}